\documentclass[10pt,a4paper]{amsart}

\usepackage{amsmath}
\usepackage{amssymb}
\usepackage{amsthm}
\usepackage{amsfonts}
\usepackage{microtype}
\usepackage{mathrsfs}
\usepackage{graphicx}
\usepackage{color}
\usepackage{latexsym}
\usepackage{rotating}
\usepackage{xspace}
\usepackage[all]{xy}
\usepackage{longtable}
\usepackage{leftidx}
\usepackage{mathtools}
\usepackage{titlesec}
\usepackage{tikz}
\usetikzlibrary{calc}
\usetikzlibrary{arrows}
\usetikzlibrary{decorations.markings}
\usepackage{geometry}

\newtheorem{theorem}{Theorem}[section]
\numberwithin{equation}{theorem}
\newtheorem{lemma}[theorem]{Lemma}
\newtheorem{proposition}[theorem]{Proposition}
\newtheorem{corollary}[theorem]{Corollary}

\theoremstyle{definition}
\newtheorem{definition}[theorem]{Definition}
\newtheorem{example}[theorem]{Example}

\newtheorem{remark}[theorem]{Remark}

\theoremstyle{conjecture}

\newcommand{\im}{\operatorname{im}}

\newcommand{\id}{\operatorname{id}}

\newcommand{\Cone}{\operatorname{Cone}}

\newcommand{\Ext}{\operatorname{Ext}}

\newcommand{\Tor}{\operatorname{Tor}}
\newcommand{\Hom}{\operatorname{Hom}}

\newcommand{\coker}{\operatorname{coker}}

\newcommand{\suchthat}{\;\ifnum\currentgrouptype=16 \middle\fi|\;}

\newenvironment{prf}[1][Proof]{\begin{proof}[\bf #1]}{\end{proof}}

\makeatletter
\newcommand{\hocolim@}[2]{%
  \vtop{\m@th\ialign{##\cr
    \hfil$#1\operator@font holim$\hfil\cr
    \noalign{\nointerlineskip\kern1.5\ex@}#2\cr
    \noalign{\nointerlineskip\kern-\ex@}\cr}}%
}
\newcommand{\hocolim}{%
  \mathop{\mathpalette\hocolim@{\rightarrowfill@\textstyle}}\nmlimits@
}
\makeatother

\makeatletter
\newcommand{\holim@}[2]{%
  \vtop{\m@th\ialign{##\cr
    \hfil$#1\operator@font holim$\hfil\cr
    \noalign{\nointerlineskip\kern1.5\ex@}#2\cr
    \noalign{\nointerlineskip\kern-\ex@}\cr}}%
}
\newcommand{\holim}{%
  \mathop{\mathpalette\holim@{\leftarrowfill@\textstyle}}\nmlimits@
}
\makeatother

\makeatletter
\def\@secnumfont{\bfseries}
\def\section{\@startsection{section}{1}%
  \z@{.7\linespacing\@plus\linespacing}{.5\linespacing}%
  {\normalfont\Large\bfseries\filcenter}}
\def\subsection{\@startsection{subsection}{2}%
  \z@{.5\linespacing\@plus.7\linespacing}{-.5em}%
  {\normalfont\large\bfseries}}
\makeatother

\begin{document}

\author[H. Faridian]{Hossein Faridian}

\title[Greenlees-May Duality in a Nutshell]
{Greenlees-May Duality in a Nutshell}

\address{H. Faridian, Department of Mathematics, Shahid Beheshti University, G.C., Evin, Tehran, Iran, Zip Code 1983963113.}
\email{h.faridian@yahoo.com}



\begin{abstract}
This expository article delves into the Greenlees-May Duality Theorem which is widely thought of as a far-reaching generalization of the Grothendieck's Local Duality Theorem. This theorem is not addressed in the literature as it merits and its proof is indeed a tangled web in a series of scattered papers. By carefully scrutinizing the requisite tools, we present a clear-cut well-documented proof of this theorem for the sake of bookkeeping.
\end{abstract}

\maketitle

\tableofcontents

\sloppy

\section{Introduction}

Throughout this note, all rings are assumed to be commutative and noetherian with identity.

In his algebraic geometry seminars of 1961-2, Grothendieck founded the theory of local cohomology as an indispensable tool in both algebraic geometry and commutative algebra. Given an ideal $\mathfrak{a}$ of $R$, the local cohomology functor $H^{i}_{\mathfrak{a}}(-)$ is defined as the $i$th right derived functor of the $\mathfrak{a}$-torsion functor $\Gamma_{\mathfrak{a}}(-) \cong \varinjlim \Hom_{R}(R/\mathfrak{a}^{n},-)$. Among a myriad of outstanding results, he proved the Local Duality Theorem.

\begin{theorem} \label{3.1.1}
Let $(R,\mathfrak{m})$ be a local ring with a dualizing module $\omega_{R}$, and $M$ a finitely generated $R$-module. Then
$$H^{i}_{\mathfrak{m}}(M) \cong \Ext^{\dim(R)-i}_{R}(M,\omega_{R})^{\vee}$$
for every $i \geq 0$.
\end{theorem}

The dual theory to local cohomology, i.e. local homology, was initiated by Matlis \cite{Mat} in 1974, and its study was continued by Simon in \cite{Si1} and \cite{Si2}. Given an ideal $\mathfrak{a}$ of $R$, the local homology functor $H^{\mathfrak{a}}_{i}(-)$ is defined as the $i$th left derived functor of the $\mathfrak{a}$-adic completion functor $\Lambda^{\mathfrak{a}}(-) \cong \varprojlim (R/\mathfrak{a}^{n}\otimes_{R}-)$.

The existence of a dualizing module in Theorem \ref{3.1.1} is rather restrictive as it forces $R$ to be Cohen-Macaulay. To proceed further and generalize Theorem \ref{3.1.1}, Greenlees and May \cite[Propositions 3.1 and 3.8]{GM}, established a spectral sequence
\begin{equation} \label{eq:3.1.1.1}
E^{2}_{p,q}= \Ext^{-p}_{R}\left(H^{q}_{\mathfrak{a}}(R),M\right) \underset {p} \Rightarrow H^{\mathfrak{a}}_{p+q}(M)
\end{equation}
for any $R$-module $M$.
One can also settle the dual spectral sequence
\begin{equation} \label{eq:3.1.1.2}
E^{2}_{p,q}= \Tor^{R}_{p}\left(H^{q}_{\mathfrak{a}}(R),M\right) \underset {p} \Rightarrow H^{p+q}_{\mathfrak{a}}(M)
\end{equation}
for any $R$-module $M$.
It is by and large more palatable to have isomorphisms rather than spectral sequences. But the problem is that the category of $R$-modules $\mathcal{M}(R)$ is not rich enough to allow for the coveted isomorphisms. We need to enlarge this category to the category of $R$-complexes $\mathcal{C}(R)$, and even enrich it further, to the derived category $\mathcal{D}(R)$. The derived category $\mathcal{D}(R)$ is privileged with extreme maturity to accommodate the sought isomorphisms. As a matter of fact, the spectral sequence \eqref{eq:3.1.1.1} turns into the isomorphism
\begin{equation} \label{eq:3.1.1.3}
{\bf R}\Hom_{R}\left({\bf R}\Gamma_{\mathfrak{a}}(R),X\right)\simeq {\bf L}\Lambda^{\mathfrak{a}}(X),
\end{equation}
and the spectral sequence \eqref{eq:3.1.1.2} turns into the isomorphism
\begin{equation} \label{eq:3.1.1.4}
{\bf R}\Gamma_{\mathfrak{a}}(R)\otimes_{R}^{\bf L}X \simeq {\bf R}\Gamma_{\mathfrak{a}}(X)
\end{equation}
in $\mathcal{D}(R)$ for any $R$-complex $X$. Patching the two isomorphisms \eqref{eq:3.1.1.3} and \eqref{eq:3.1.1.4} together, we are blessed with the celebrated Greenlees-May Duality.
\begin{theorem} \label{3.1.2}
Let $\mathfrak{a}$ be an ideal of $R$, and $X,Y \in \mathcal{D}(R)$. Then there is a natural isomorphism
$${\bf R}\Hom_{R}\left({\bf R}\Gamma_{\mathfrak{a}}(X),Y\right) \simeq {\bf R}\Hom_{R}\left(X,{\bf L}\Lambda^{\mathfrak{a}}(Y)\right)$$
in $\mathcal{D}(R)$.
\end{theorem}
This was first proved by Alonso Tarr\'{i}o, Jerem\'{i}as L\'{o}pez and Lipman in \cite{AJL}. Theorem \ref{3.1.2} is a far-reaching generalization of Theorem \ref{3.1.1} and indeed extends it to its full generality. This theorem also demonstrates perfectly some sort of adjointness between derived local cohomology and homology.

Despite its incontrovertible impact on the theory of derived local homology and cohomology, we regretfully notice that there is no comprehensive and accessible treatment of the Greenlees-May Duality in the literature. There are some papers that touch on the subject, each from a different perspective, but none of them present a clear-cut and thorough proof that is fairly readable for non-experts; see for example \cite{GM}, \cite{AJL}, \cite{PSY}, and \cite{Sc}. To remedy this defect, we commence on probing this theorem by providing the prerequisites from scratch and build upon a well-documented rigorous proof which is basically presented in layman's terms. In the course of our proof, some arguments are familiar while some others are novel. However, all the details are fully worked out so as to set forth a satisfactory exposition of the subject. We finally depict the highly non-trivial fact that the Greenlees-May Duality generalizes the Local Duality in simple and traceable steps.

\section{Module Prerequisites}

In this section, we embark on providing the requisite tools on modules which are to be recruited in Section 4.

First we recall the notion of a $\delta$-functor which will be used as a powerful tool to establish natural isomorphisms.

\begin{definition} \label{3.2.1}
Let $R$ and $S$ be two rings. Then:
\begin{enumerate}
\item[(i)] A \textit{homological covariant $\delta$-functor} \index{homological covariant $\delta$-functor} is a sequence $\left(\mathcal{F}_{i}:\mathcal{M}(R)\rightarrow \mathcal{M}(S)\right)_{i \geq 0}$ of additive covariant functors with the property that every short exact sequence
$$0 \rightarrow M' \rightarrow M \rightarrow M'' \rightarrow 0$$
of $R$-modules gives rise to a long exact sequence
$$\cdots \rightarrow \mathcal{F}_{2}(M'') \xrightarrow{\delta_{2}} \mathcal{F}_{1}(M') \rightarrow \mathcal{F}_{1}(M) \rightarrow \mathcal{F}_{1}(M'') \xrightarrow{\delta_{1}} \mathcal{F}_{0}(M') \rightarrow \mathcal{F}_{0}(M) \rightarrow \mathcal{F}_{0}(M'') \rightarrow 0$$
of $S$-modules, such that the connecting morphisms $\delta_{i}$'s are natural in the sense that any commutative diagram
\[\begin{tikzpicture}[every node/.style={midway},]
  \matrix[column sep={2.5em}, row sep={2.5em}]
  {\node(1) {$0$}; & \node(2) {$M'$}; & \node(3) {$M$}; & \node(4) {$M''$}; & \node(5) {$0$};\\
  \node(6) {$0$}; & \node(7) {$N'$}; & \node(8) {$N$}; & \node(9) {$N''$}; & \node(10) {$0$};\\};
  \draw[decoration={markings,mark=at position 1 with {\arrow[scale=1.5]{>}}},postaction={decorate},shorten >=0.5pt] (4) -- (9) node[anchor=east] {};
  \draw[decoration={markings,mark=at position 1 with {\arrow[scale=1.5]{>}}},postaction={decorate},shorten >=0.5pt] (3) -- (8) node[anchor=east] {};
  \draw[decoration={markings,mark=at position 1 with {\arrow[scale=1.5]{>}}},postaction={decorate},shorten >=0.5pt] (2) -- (7) node[anchor=east] {};
  \draw[decoration={markings,mark=at position 1 with {\arrow[scale=1.5]{>}}},postaction={decorate},shorten >=0.5pt] (1) -- (2) node[anchor=south] {};
  \draw[decoration={markings,mark=at position 1 with {\arrow[scale=1.5]{>}}},postaction={decorate},shorten >=0.5pt] (2) -- (3) node[anchor=south] {};
  \draw[decoration={markings,mark=at position 1 with {\arrow[scale=1.5]{>}}},postaction={decorate},shorten >=0.5pt] (3) -- (4) node[anchor=south] {};
  \draw[decoration={markings,mark=at position 1 with {\arrow[scale=1.5]{>}}},postaction={decorate},shorten >=0.5pt] (4) -- (5) node[anchor=south] {};
  \draw[decoration={markings,mark=at position 1 with {\arrow[scale=1.5]{>}}},postaction={decorate},shorten >=0.5pt] (6) -- (7) node[anchor=south] {};
  \draw[decoration={markings,mark=at position 1 with {\arrow[scale=1.5]{>}}},postaction={decorate},shorten >=0.5pt] (7) -- (8) node[anchor=south] {};
  \draw[decoration={markings,mark=at position 1 with {\arrow[scale=1.5]{>}}},postaction={decorate},shorten >=0.5pt] (8) -- (9) node[anchor=south] {};
  \draw[decoration={markings,mark=at position 1 with {\arrow[scale=1.5]{>}}},postaction={decorate},shorten >=0.5pt] (9) -- (10) node[anchor=south] {};
\end{tikzpicture}\]
of $R$-modules with exact rows induces a commutative diagram
\[\begin{tikzpicture}[every node/.style={midway},]
  \matrix[column sep={1em}, row sep={2.5em}]
  {\node(1) {$\cdots$}; & \node(2) {$\mathcal{F}_{2}(M'')$}; & \node(3) {$\mathcal{F}_{1}(M')$}; & \node(4) {$\mathcal{F}_{1}(M)$}; & \node(5) {$\mathcal{F}_{1}(M'')$}; & \node(6) {$\mathcal{F}_{0}(M')$}; & \node(7) {$\mathcal{F}_{0}(M)$}; & \node(8) {$\mathcal{F}_{0}(M'')$}; & \node(9) {$0$};\\
  \node(10) {$\cdots$}; & \node(11) {$\mathcal{F}_{2}(N'')$}; & \node(12) {$\mathcal{F}_{1}(N')$}; & \node(13) {$\mathcal{F}_{1}(N)$}; & \node(14) {$\mathcal{F}_{1}(N'')$}; & \node(15) {$\mathcal{F}_{0}(N')$}; & \node(16) {$\mathcal{F}_{0}(N)$}; & \node(17) {$\mathcal{F}_{0}(N'')$}; & \node(18) {$0$};\\};
  \draw[decoration={markings,mark=at position 1 with {\arrow[scale=1.5]{>}}},postaction={decorate},shorten >=0.5pt] (1) -- (2) node[anchor=east] {};
  \draw[decoration={markings,mark=at position 1 with {\arrow[scale=1.5]{>}}},postaction={decorate},shorten >=0.5pt] (2) -- (3) node[anchor=south] {$\delta_{2}$};
  \draw[decoration={markings,mark=at position 1 with {\arrow[scale=1.5]{>}}},postaction={decorate},shorten >=0.5pt] (3) -- (4) node[anchor=east] {};
  \draw[decoration={markings,mark=at position 1 with {\arrow[scale=1.5]{>}}},postaction={decorate},shorten >=0.5pt] (4) -- (5) node[anchor=east] {};
  \draw[decoration={markings,mark=at position 1 with {\arrow[scale=1.5]{>}}},postaction={decorate},shorten >=0.5pt] (5) -- (6) node[anchor=south] {$\delta_{1}$};
  \draw[decoration={markings,mark=at position 1 with {\arrow[scale=1.5]{>}}},postaction={decorate},shorten >=0.5pt] (6) -- (7) node[anchor=east] {};
  \draw[decoration={markings,mark=at position 1 with {\arrow[scale=1.5]{>}}},postaction={decorate},shorten >=0.5pt] (7) -- (8) node[anchor=east] {};
  \draw[decoration={markings,mark=at position 1 with {\arrow[scale=1.5]{>}}},postaction={decorate},shorten >=0.5pt] (8) -- (9) node[anchor=east] {};
  \draw[decoration={markings,mark=at position 1 with {\arrow[scale=1.5]{>}}},postaction={decorate},shorten >=0.5pt] (10) -- (11) node[anchor=east] {};
  \draw[decoration={markings,mark=at position 1 with {\arrow[scale=1.5]{>}}},postaction={decorate},shorten >=0.5pt] (11) -- (12) node[anchor=south] {$\Delta_{2}$};
  \draw[decoration={markings,mark=at position 1 with {\arrow[scale=1.5]{>}}},postaction={decorate},shorten >=0.5pt] (12) -- (13) node[anchor=east] {};
  \draw[decoration={markings,mark=at position 1 with {\arrow[scale=1.5]{>}}},postaction={decorate},shorten >=0.5pt] (13) -- (14) node[anchor=east] {};
  \draw[decoration={markings,mark=at position 1 with {\arrow[scale=1.5]{>}}},postaction={decorate},shorten >=0.5pt] (14) -- (15) node[anchor=south] {$\Delta_{1}$};
  \draw[decoration={markings,mark=at position 1 with {\arrow[scale=1.5]{>}}},postaction={decorate},shorten >=0.5pt] (15) -- (16) node[anchor=east] {};
  \draw[decoration={markings,mark=at position 1 with {\arrow[scale=1.5]{>}}},postaction={decorate},shorten >=0.5pt] (16) -- (17) node[anchor=east] {};
  \draw[decoration={markings,mark=at position 1 with {\arrow[scale=1.5]{>}}},postaction={decorate},shorten >=0.5pt] (17) -- (18) node[anchor=east] {};
  \draw[decoration={markings,mark=at position 1 with {\arrow[scale=1.5]{>}}},postaction={decorate},shorten >=0.5pt] (2) -- (11) node[anchor=east] {};
  \draw[decoration={markings,mark=at position 1 with {\arrow[scale=1.5]{>}}},postaction={decorate},shorten >=0.5pt] (3) -- (12) node[anchor=east] {};
  \draw[decoration={markings,mark=at position 1 with {\arrow[scale=1.5]{>}}},postaction={decorate},shorten >=0.5pt] (4) -- (13) node[anchor=east] {};
  \draw[decoration={markings,mark=at position 1 with {\arrow[scale=1.5]{>}}},postaction={decorate},shorten >=0.5pt] (5) -- (14) node[anchor=east] {};
  \draw[decoration={markings,mark=at position 1 with {\arrow[scale=1.5]{>}}},postaction={decorate},shorten >=0.5pt] (6) -- (15) node[anchor=east] {};
  \draw[decoration={markings,mark=at position 1 with {\arrow[scale=1.5]{>}}},postaction={decorate},shorten >=0.5pt] (7) -- (16) node[anchor=east] {};
  \draw[decoration={markings,mark=at position 1 with {\arrow[scale=1.5]{>}}},postaction={decorate},shorten >=0.5pt] (8) -- (17) node[anchor=east] {};
\end{tikzpicture}\]
of $S$-modules with exact rows.
\item[(ii)] A \textit{cohomological covariant $\delta$-functor} \index{cohomological covariant $\delta$-functor} is a sequence $\left(\mathcal{F}^{i}:\mathcal{M}(R)\rightarrow \mathcal{M}(S)\right)_{i \geq 0}$ of additive covariant functors with the property that every short exact sequence
$$0 \rightarrow M' \rightarrow M \rightarrow M'' \rightarrow 0$$
of $R$-modules gives rise to a long exact sequence
$$0 \rightarrow \mathcal{F}^{0}(M') \rightarrow \mathcal{F}^{0}(M) \rightarrow \mathcal{F}^{0}(M'') \xrightarrow{\delta^{0}} \mathcal{F}^{1}(M') \rightarrow \mathcal{F}^{1}(M) \rightarrow \mathcal{F}^{1}(M'') \xrightarrow{\delta^{1}} \mathcal{F}^{2}(M') \rightarrow \cdots$$
of $S$-modules, such that the connecting morphisms $\delta^{i}$'s are natural in the sense that any commutative diagram
\[\begin{tikzpicture}[every node/.style={midway},]
  \matrix[column sep={2.5em}, row sep={2.5em}]
  {\node(1) {$0$}; & \node(2) {$M'$}; & \node(3) {$M$}; & \node(4) {$M''$}; & \node(5) {$0$};\\
  \node(6) {$0$}; & \node(7) {$N'$}; & \node(8) {$N$}; & \node(9) {$N''$}; & \node(10) {$0$};\\};
  \draw[decoration={markings,mark=at position 1 with {\arrow[scale=1.5]{>}}},postaction={decorate},shorten >=0.5pt] (4) -- (9) node[anchor=east] {};
  \draw[decoration={markings,mark=at position 1 with {\arrow[scale=1.5]{>}}},postaction={decorate},shorten >=0.5pt] (3) -- (8) node[anchor=east] {};
  \draw[decoration={markings,mark=at position 1 with {\arrow[scale=1.5]{>}}},postaction={decorate},shorten >=0.5pt] (2) -- (7) node[anchor=east] {};
  \draw[decoration={markings,mark=at position 1 with {\arrow[scale=1.5]{>}}},postaction={decorate},shorten >=0.5pt] (1) -- (2) node[anchor=south] {};
  \draw[decoration={markings,mark=at position 1 with {\arrow[scale=1.5]{>}}},postaction={decorate},shorten >=0.5pt] (2) -- (3) node[anchor=south] {};
  \draw[decoration={markings,mark=at position 1 with {\arrow[scale=1.5]{>}}},postaction={decorate},shorten >=0.5pt] (3) -- (4) node[anchor=south] {};
  \draw[decoration={markings,mark=at position 1 with {\arrow[scale=1.5]{>}}},postaction={decorate},shorten >=0.5pt] (4) -- (5) node[anchor=south] {};
  \draw[decoration={markings,mark=at position 1 with {\arrow[scale=1.5]{>}}},postaction={decorate},shorten >=0.5pt] (6) -- (7) node[anchor=south] {};
  \draw[decoration={markings,mark=at position 1 with {\arrow[scale=1.5]{>}}},postaction={decorate},shorten >=0.5pt] (7) -- (8) node[anchor=south] {};
  \draw[decoration={markings,mark=at position 1 with {\arrow[scale=1.5]{>}}},postaction={decorate},shorten >=0.5pt] (8) -- (9) node[anchor=south] {};
  \draw[decoration={markings,mark=at position 1 with {\arrow[scale=1.5]{>}}},postaction={decorate},shorten >=0.5pt] (9) -- (10) node[anchor=south] {};
\end{tikzpicture}\]
of $R$-modules with exact rows induces a commutative diagram
\[\begin{tikzpicture}[every node/.style={midway},]
  \matrix[column sep={1em}, row sep={2.5em}]
  {\node(1) {$0$}; & \node(2) {$\mathcal{F}^{0}(M')$}; & \node(3) {$\mathcal{F}^{0}(M)$}; & \node(4) {$\mathcal{F}^{0}(M'')$}; & \node(5) {$\mathcal{F}^{1}(M')$}; & \node(6) {$\mathcal{F}^{1}(M)$}; & \node(7) {$\mathcal{F}^{1}(M'')$}; & \node(8) {$\mathcal{F}^{2}(M')$}; & \node(9) {$\cdots$};\\
  \node(10) {$0$}; & \node(11) {$\mathcal{F}^{0}(N')$}; & \node(12) {$\mathcal{F}^{0}(N)$}; & \node(13) {$\mathcal{F}^{0}(N'')$}; & \node(14) {$\mathcal{F}^{1}(N')$}; & \node(15) {$\mathcal{F}^{1}(N)$}; & \node(16) {$\mathcal{F}^{1}(N'')$}; & \node(17) {$\mathcal{F}^{2}(N')$}; & \node(18) {$\cdots$};\\};
  \draw[decoration={markings,mark=at position 1 with {\arrow[scale=1.5]{>}}},postaction={decorate},shorten >=0.5pt] (1) -- (2) node[anchor=east] {};
  \draw[decoration={markings,mark=at position 1 with {\arrow[scale=1.5]{>}}},postaction={decorate},shorten >=0.5pt] (2) -- (3) node[anchor=south] {};
  \draw[decoration={markings,mark=at position 1 with {\arrow[scale=1.5]{>}}},postaction={decorate},shorten >=0.5pt] (3) -- (4) node[anchor=east] {};
  \draw[decoration={markings,mark=at position 1 with {\arrow[scale=1.5]{>}}},postaction={decorate},shorten >=0.5pt] (4) -- (5) node[anchor=south] {$\delta^{0}$};
  \draw[decoration={markings,mark=at position 1 with {\arrow[scale=1.5]{>}}},postaction={decorate},shorten >=0.5pt] (5) -- (6) node[anchor=south] {};
  \draw[decoration={markings,mark=at position 1 with {\arrow[scale=1.5]{>}}},postaction={decorate},shorten >=0.5pt] (6) -- (7) node[anchor=east] {};
  \draw[decoration={markings,mark=at position 1 with {\arrow[scale=1.5]{>}}},postaction={decorate},shorten >=0.5pt] (7) -- (8) node[anchor=south] {$\delta^{1}$};
  \draw[decoration={markings,mark=at position 1 with {\arrow[scale=1.5]{>}}},postaction={decorate},shorten >=0.5pt] (8) -- (9) node[anchor=east] {};
  \draw[decoration={markings,mark=at position 1 with {\arrow[scale=1.5]{>}}},postaction={decorate},shorten >=0.5pt] (10) -- (11) node[anchor=east] {};
  \draw[decoration={markings,mark=at position 1 with {\arrow[scale=1.5]{>}}},postaction={decorate},shorten >=0.5pt] (11) -- (12) node[anchor=south] {};
  \draw[decoration={markings,mark=at position 1 with {\arrow[scale=1.5]{>}}},postaction={decorate},shorten >=0.5pt] (12) -- (13) node[anchor=east] {};
  \draw[decoration={markings,mark=at position 1 with {\arrow[scale=1.5]{>}}},postaction={decorate},shorten >=0.5pt] (13) -- (14) node[anchor=south] {$\Delta^{0}$};
  \draw[decoration={markings,mark=at position 1 with {\arrow[scale=1.5]{>}}},postaction={decorate},shorten >=0.5pt] (14) -- (15) node[anchor=south] {};
  \draw[decoration={markings,mark=at position 1 with {\arrow[scale=1.5]{>}}},postaction={decorate},shorten >=0.5pt] (15) -- (16) node[anchor=east] {};
  \draw[decoration={markings,mark=at position 1 with {\arrow[scale=1.5]{>}}},postaction={decorate},shorten >=0.5pt] (16) -- (17) node[anchor=south] {$\Delta^{1}$};
  \draw[decoration={markings,mark=at position 1 with {\arrow[scale=1.5]{>}}},postaction={decorate},shorten >=0.5pt] (17) -- (18) node[anchor=east] {};
  \draw[decoration={markings,mark=at position 1 with {\arrow[scale=1.5]{>}}},postaction={decorate},shorten >=0.5pt] (2) -- (11) node[anchor=east] {};
  \draw[decoration={markings,mark=at position 1 with {\arrow[scale=1.5]{>}}},postaction={decorate},shorten >=0.5pt] (3) -- (12) node[anchor=east] {};
  \draw[decoration={markings,mark=at position 1 with {\arrow[scale=1.5]{>}}},postaction={decorate},shorten >=0.5pt] (4) -- (13) node[anchor=east] {};
  \draw[decoration={markings,mark=at position 1 with {\arrow[scale=1.5]{>}}},postaction={decorate},shorten >=0.5pt] (5) -- (14) node[anchor=east] {};
  \draw[decoration={markings,mark=at position 1 with {\arrow[scale=1.5]{>}}},postaction={decorate},shorten >=0.5pt] (6) -- (15) node[anchor=east] {};
  \draw[decoration={markings,mark=at position 1 with {\arrow[scale=1.5]{>}}},postaction={decorate},shorten >=0.5pt] (7) -- (16) node[anchor=east] {};
  \draw[decoration={markings,mark=at position 1 with {\arrow[scale=1.5]{>}}},postaction={decorate},shorten >=0.5pt] (8) -- (17) node[anchor=east] {};
\end{tikzpicture}\]
of $S$-modules with exact rows.
\end{enumerate}
\end{definition}

\begin{example} \label{3.2.2}
Let $R$ and $S$ be two rings, and $\mathcal{F}:\mathcal{M}(R) \rightarrow \mathcal{M}(S)$ an additive covariant functor. Then $\left(L_{i}\mathcal{F}:\mathcal{M}(R) \rightarrow \mathcal{M}(S)\right)_{i \geq 0}$ is a homological covariant $\delta$-functor and $\left(R^{i}\mathcal{F}:\mathcal{M}(R) \rightarrow \mathcal{M}(S)\right)_{i \geq 0}$ is a cohomological covariant $\delta$-functor.
\end{example}

\begin{definition} \label{3.2.3}
Let $R$ and $S$ be two rings. Then:
\begin{enumerate}
\item[(i)] A \textit{morphism} \index{morphism of homological covariant $\delta$-functors}
$$\tau: \left(\mathcal{F}_{i}:\mathcal{M}(R) \rightarrow \mathcal{M}(S)\right)_{i \geq 0} \rightarrow \left(\mathcal{G}_{i}:\mathcal{M}(R) \rightarrow \mathcal{M}(S)\right)_{i \geq 0}$$
of homological covariant $\delta$-functors is a sequence $\tau= \left(\tau_{i}:\mathcal{F}_{i} \rightarrow \mathcal{G}_{i}\right)_{i \geq 0}$ of natural transformations of functors, such that any short exact sequence
$$0 \rightarrow M' \rightarrow M \rightarrow M'' \rightarrow 0$$
of $R$-modules induces a commutative diagram
\[\begin{tikzpicture}[every node/.style={midway},]
  \matrix[column sep={1em}, row sep={2.5em}]
  {\node(1) {$\cdots$}; & \node(2) {$\mathcal{F}_{2}(M'')$}; & \node(3) {$\mathcal{F}_{1}(M')$}; & \node(4) {$\mathcal{F}_{1}(M)$}; & \node(5) {$\mathcal{F}_{1}(M'')$}; & \node(6) {$\mathcal{F}_{0}(M')$}; & \node(7) {$\mathcal{F}_{0}(M)$}; & \node(8) {$\mathcal{F}_{0}(M'')$}; & \node(9) {$0$};\\
  \node(10) {$\cdots$}; & \node(11) {$\mathcal{G}_{2}(M'')$}; & \node(12) {$\mathcal{G}_{1}(M')$}; & \node(13) {$\mathcal{G}_{1}(M)$}; & \node(14) {$\mathcal{G}_{1}(M'')$}; & \node(15) {$\mathcal{G}_{0}(M')$}; & \node(16) {$\mathcal{G}_{0}(M)$}; & \node(17) {$\mathcal{G}_{0}(M'')$}; & \node(18) {$0$};\\};
  \draw[decoration={markings,mark=at position 1 with {\arrow[scale=1.5]{>}}},postaction={decorate},shorten >=0.5pt] (1) -- (2) node[anchor=east] {};
  \draw[decoration={markings,mark=at position 1 with {\arrow[scale=1.5]{>}}},postaction={decorate},shorten >=0.5pt] (2) -- (3) node[anchor=south] {$\delta_{2}$};
  \draw[decoration={markings,mark=at position 1 with {\arrow[scale=1.5]{>}}},postaction={decorate},shorten >=0.5pt] (3) -- (4) node[anchor=east] {};
  \draw[decoration={markings,mark=at position 1 with {\arrow[scale=1.5]{>}}},postaction={decorate},shorten >=0.5pt] (4) -- (5) node[anchor=east] {};
  \draw[decoration={markings,mark=at position 1 with {\arrow[scale=1.5]{>}}},postaction={decorate},shorten >=0.5pt] (5) -- (6) node[anchor=south] {$\delta_{1}$};
  \draw[decoration={markings,mark=at position 1 with {\arrow[scale=1.5]{>}}},postaction={decorate},shorten >=0.5pt] (6) -- (7) node[anchor=east] {};
  \draw[decoration={markings,mark=at position 1 with {\arrow[scale=1.5]{>}}},postaction={decorate},shorten >=0.5pt] (7) -- (8) node[anchor=east] {};
  \draw[decoration={markings,mark=at position 1 with {\arrow[scale=1.5]{>}}},postaction={decorate},shorten >=0.5pt] (8) -- (9) node[anchor=east] {};
  \draw[decoration={markings,mark=at position 1 with {\arrow[scale=1.5]{>}}},postaction={decorate},shorten >=0.5pt] (10) -- (11) node[anchor=east] {};
  \draw[decoration={markings,mark=at position 1 with {\arrow[scale=1.5]{>}}},postaction={decorate},shorten >=0.5pt] (11) -- (12) node[anchor=south] {$\Delta_{2}$};
  \draw[decoration={markings,mark=at position 1 with {\arrow[scale=1.5]{>}}},postaction={decorate},shorten >=0.5pt] (12) -- (13) node[anchor=east] {};
  \draw[decoration={markings,mark=at position 1 with {\arrow[scale=1.5]{>}}},postaction={decorate},shorten >=0.5pt] (13) -- (14) node[anchor=east] {};
  \draw[decoration={markings,mark=at position 1 with {\arrow[scale=1.5]{>}}},postaction={decorate},shorten >=0.5pt] (14) -- (15) node[anchor=south] {$\Delta_{1}$};
  \draw[decoration={markings,mark=at position 1 with {\arrow[scale=1.5]{>}}},postaction={decorate},shorten >=0.5pt] (15) -- (16) node[anchor=east] {};
  \draw[decoration={markings,mark=at position 1 with {\arrow[scale=1.5]{>}}},postaction={decorate},shorten >=0.5pt] (16) -- (17) node[anchor=east] {};
  \draw[decoration={markings,mark=at position 1 with {\arrow[scale=1.5]{>}}},postaction={decorate},shorten >=0.5pt] (17) -- (18) node[anchor=east] {};
  \draw[decoration={markings,mark=at position 1 with {\arrow[scale=1.5]{>}}},postaction={decorate},shorten >=0.5pt] (2) -- (11) node[anchor=west] {$\tau_{2}(M'')$};
  \draw[decoration={markings,mark=at position 1 with {\arrow[scale=1.5]{>}}},postaction={decorate},shorten >=0.5pt] (3) -- (12) node[anchor=west] {$\tau_{1}(M')$};
  \draw[decoration={markings,mark=at position 1 with {\arrow[scale=1.5]{>}}},postaction={decorate},shorten >=0.5pt] (4) -- (13) node[anchor=west] {$\tau_{1}(M)$};
  \draw[decoration={markings,mark=at position 1 with {\arrow[scale=1.5]{>}}},postaction={decorate},shorten >=0.5pt] (5) -- (14) node[anchor=west] {$\tau_{1}(M'')$};
  \draw[decoration={markings,mark=at position 1 with {\arrow[scale=1.5]{>}}},postaction={decorate},shorten >=0.5pt] (6) -- (15) node[anchor=west] {$\tau_{0}(M')$};
  \draw[decoration={markings,mark=at position 1 with {\arrow[scale=1.5]{>}}},postaction={decorate},shorten >=0.5pt] (7) -- (16) node[anchor=west] {$\tau_{0}(M)$};
  \draw[decoration={markings,mark=at position 1 with {\arrow[scale=1.5]{>}}},postaction={decorate},shorten >=0.5pt] (8) -- (17) node[anchor=west] {$\tau_{0}(M'')$};
\end{tikzpicture}\]
of $S$-modules with exact rows. If in particular, $\tau_{i}$ is an isomorphism for every $i \geq 0$, then $\tau$ is called an \textit{isomorphism} \index{isomorphism of homological covariant $\delta$-functors} of $\delta$-functors.
\item[(ii)] A \textit{morphism} \index{morphism of cohomological covariant $\delta$-functors}
$$\tau: \left(\mathcal{F}^{i}:\mathcal{M}(R) \rightarrow \mathcal{M}(S)\right)_{i \geq 0} \rightarrow \left(\mathcal{G}^{i}:\mathcal{M}(R) \rightarrow \mathcal{M}(S)\right)_{i \geq 0}$$
of cohomological covariant $\delta$-functors is a sequence $\tau= \left(\tau^{i}:\mathcal{F}^{i} \rightarrow \mathcal{G}^{i}\right)_{i \geq 0}$ of natural transformations of functors, such that any short exact sequence
$$0 \rightarrow M' \rightarrow M \rightarrow M'' \rightarrow 0$$
of $R$-modules induces a commutative diagram
\[\begin{tikzpicture}[every node/.style={midway},]
  \matrix[column sep={1em}, row sep={2.5em}]
  {\node(1) {$0$}; & \node(2) {$\mathcal{F}^{0}(M')$}; & \node(3) {$\mathcal{F}^{0}(M)$}; & \node(4) {$\mathcal{F}^{0}(M'')$}; & \node(5) {$\mathcal{F}^{1}(M')$}; & \node(6) {$\mathcal{F}^{1}(M)$}; & \node(7) {$\mathcal{F}^{1}(M'')$}; & \node(8) {$\mathcal{F}^{2}(M')$}; & \node(9) {$\cdots$};\\
  \node(10) {$0$}; & \node(11) {$\mathcal{G}^{0}(M')$}; & \node(12) {$\mathcal{G}^{0}(M)$}; & \node(13) {$\mathcal{G}^{0}(M'')$}; & \node(14) {$\mathcal{G}^{1}(M')$}; & \node(15) {$\mathcal{G}^{1}(M)$}; & \node(16) {$\mathcal{G}^{1}(M'')$}; & \node(17) {$\mathcal{G}^{2}(M')$}; & \node(18) {$\cdots$};\\};
  \draw[decoration={markings,mark=at position 1 with {\arrow[scale=1.5]{>}}},postaction={decorate},shorten >=0.5pt] (1) -- (2) node[anchor=east] {};
  \draw[decoration={markings,mark=at position 1 with {\arrow[scale=1.5]{>}}},postaction={decorate},shorten >=0.5pt] (2) -- (3) node[anchor=south] {};
  \draw[decoration={markings,mark=at position 1 with {\arrow[scale=1.5]{>}}},postaction={decorate},shorten >=0.5pt] (3) -- (4) node[anchor=east] {};
  \draw[decoration={markings,mark=at position 1 with {\arrow[scale=1.5]{>}}},postaction={decorate},shorten >=0.5pt] (4) -- (5) node[anchor=south] {$\delta^{0}$};
  \draw[decoration={markings,mark=at position 1 with {\arrow[scale=1.5]{>}}},postaction={decorate},shorten >=0.5pt] (5) -- (6) node[anchor=south] {};
  \draw[decoration={markings,mark=at position 1 with {\arrow[scale=1.5]{>}}},postaction={decorate},shorten >=0.5pt] (6) -- (7) node[anchor=east] {};
  \draw[decoration={markings,mark=at position 1 with {\arrow[scale=1.5]{>}}},postaction={decorate},shorten >=0.5pt] (7) -- (8) node[anchor=south] {$\delta^{1}$};
  \draw[decoration={markings,mark=at position 1 with {\arrow[scale=1.5]{>}}},postaction={decorate},shorten >=0.5pt] (8) -- (9) node[anchor=east] {};
  \draw[decoration={markings,mark=at position 1 with {\arrow[scale=1.5]{>}}},postaction={decorate},shorten >=0.5pt] (10) -- (11) node[anchor=east] {};
  \draw[decoration={markings,mark=at position 1 with {\arrow[scale=1.5]{>}}},postaction={decorate},shorten >=0.5pt] (11) -- (12) node[anchor=south] {};
  \draw[decoration={markings,mark=at position 1 with {\arrow[scale=1.5]{>}}},postaction={decorate},shorten >=0.5pt] (12) -- (13) node[anchor=east] {};
  \draw[decoration={markings,mark=at position 1 with {\arrow[scale=1.5]{>}}},postaction={decorate},shorten >=0.5pt] (13) -- (14) node[anchor=south] {$\Delta^{0}$};
  \draw[decoration={markings,mark=at position 1 with {\arrow[scale=1.5]{>}}},postaction={decorate},shorten >=0.5pt] (14) -- (15) node[anchor=south] {};
  \draw[decoration={markings,mark=at position 1 with {\arrow[scale=1.5]{>}}},postaction={decorate},shorten >=0.5pt] (15) -- (16) node[anchor=east] {};
  \draw[decoration={markings,mark=at position 1 with {\arrow[scale=1.5]{>}}},postaction={decorate},shorten >=0.5pt] (16) -- (17) node[anchor=south] {$\Delta^{1}$};
  \draw[decoration={markings,mark=at position 1 with {\arrow[scale=1.5]{>}}},postaction={decorate},shorten >=0.5pt] (17) -- (18) node[anchor=east] {};
  \draw[decoration={markings,mark=at position 1 with {\arrow[scale=1.5]{>}}},postaction={decorate},shorten >=0.5pt] (2) -- (11) node[anchor=west] {$\tau^{0}(M')$};
  \draw[decoration={markings,mark=at position 1 with {\arrow[scale=1.5]{>}}},postaction={decorate},shorten >=0.5pt] (3) -- (12) node[anchor=west] {$\tau^{0}(M)$};
  \draw[decoration={markings,mark=at position 1 with {\arrow[scale=1.5]{>}}},postaction={decorate},shorten >=0.5pt] (4) -- (13) node[anchor=west] {$\tau^{0}(M'')$};
  \draw[decoration={markings,mark=at position 1 with {\arrow[scale=1.5]{>}}},postaction={decorate},shorten >=0.5pt] (5) -- (14) node[anchor=west] {$\tau^{1}(M')$};
  \draw[decoration={markings,mark=at position 1 with {\arrow[scale=1.5]{>}}},postaction={decorate},shorten >=0.5pt] (6) -- (15) node[anchor=west] {$\tau^{1}(M)$};
  \draw[decoration={markings,mark=at position 1 with {\arrow[scale=1.5]{>}}},postaction={decorate},shorten >=0.5pt] (7) -- (16) node[anchor=west] {$\tau^{1}(M'')$};
  \draw[decoration={markings,mark=at position 1 with {\arrow[scale=1.5]{>}}},postaction={decorate},shorten >=0.5pt] (8) -- (17) node[anchor=west] {$\tau^{2}(M')$};
\end{tikzpicture}\]
of $S$-modules with exact rows. If in particular, $\tau^{i}$ is an isomorphism for every $i \geq 0$, then $\tau$ is called an \textit{isomorphism} \index{isomorphism of homological covariant $\delta$-functors} of $\delta$-functors.
\end{enumerate}
\end{definition}

The following remarkable theorem due to Grothendieck provides hands-on conditions that ascertain the existence of isomorphisms between $\delta$-functors.

\begin{theorem} \label{3.2.4}
Let $R$ and $S$ be two rings. Then the following assertions hold:
\begin{enumerate}
\item[(i)] Assume that $\left(\mathcal{F}_{i}:\mathcal{M}(R) \rightarrow \mathcal{M}(S)\right)_{i \geq 0}$ and $\left(\mathcal{G}_{i}:\mathcal{M}(R) \rightarrow \mathcal{M}(S)\right)_{i \geq 0}$ are two homological covariant $\delta$-functors such that $\mathcal{F}_{i}(F)=0=\mathcal{G}_{i}(F)$ for every free $R$-module $F$ and every $i \geq 1$. If there is a natural transformation $\eta: \mathcal{F}_{0} \rightarrow \mathcal{G}_{0}$ of functors which is an isomorphism on free $R$-modules, then there is a unique isomorphism $\tau:(\mathcal{F}_{i})_{i \geq 0} \rightarrow (\mathcal{G}_{i})_{i \geq 0}$ of $\delta$-functors such that $\tau_{0}=\eta$.
\item[(ii)] Assume that $\left(\mathcal{F}^{i}:\mathcal{M}(R) \rightarrow \mathcal{M}(S)\right)_{i \geq 0}$ and $\left(\mathcal{G}^{i}:\mathcal{M}(R) \rightarrow \mathcal{M}(S)\right)_{i \geq 0}$ are two cohomological covariant $\delta$-functors such that $\mathcal{F}^{i}(I)=0=\mathcal{G}^{i}(I)$ for every injective $R$-module $I$ and every $i \geq 1$. If there is a natural transformation $\eta: \mathcal{F}^{0} \rightarrow \mathcal{G}^{0}$ of functors which is an isomorphism on injective $R$-modules, then there is a unique isomorphism $\tau:(\mathcal{F}^{i})_{i \geq 0} \rightarrow (\mathcal{G}^{i})_{i \geq 0}$ of $\delta$-functors such that $\tau^{0}=\eta$.
\end{enumerate}
\end{theorem}

\begin{prf}
The proof is standard and can be found in almost every book on homological algebra. For example, see \cite[Corollaries 6.34 and 6.49]{Ro}. One should note that the above version is somewhat stronger than what is normally recorded in the books. However, the same proof can be modified in a suitable way to imply the above version.
\end{prf}

The following corollary sets forth a special case of Theorem \ref{3.2.4} which frequently occurs in practice.

\begin{corollary} \label{3.2.5}
Let $R$ and $S$ be two rings. Then the following assertions hold:
\begin{enumerate}
\item[(i)] Assume that $\mathcal{F}:\mathcal{M}(R) \rightarrow \mathcal{M}(S)$ is an additive covariant functor, and $\left(\mathcal{F}_{i}:\mathcal{M}(R) \rightarrow \mathcal{M}(S)\right)_{i \geq 0}$ is a homological covariant $\delta$-functor such that $\mathcal{F}_{i}(F)=0$ for every free $R$-module $F$ and every $i \geq 1$. If there is a natural transformation $\eta: L_{0}\mathcal{F} \rightarrow \mathcal{F}_{0}$ of functors which is an isomorphism on free $R$-modules, then there is a unique isomorphism $\tau:(L_{i}\mathcal{F})_{i \geq 0} \rightarrow (\mathcal{F}_{i})_{i \geq 0}$ of $\delta$-functors such that $\tau_{0}=\eta$.
\item[(ii)] Assume that $\mathcal{F}:\mathcal{M}(R) \rightarrow \mathcal{M}(S)$ is an additive covariant functor, and $\left(\mathcal{F}^{i}:\mathcal{M}(R) \rightarrow \mathcal{M}(S)\right)_{i \geq 0}$ is a cohomological covariant $\delta$-functor such that $\mathcal{F}^{i}(I)=0$ for every injective $R$-module $I$ and every $i \geq 1$. If there is a natural transformation $\eta: \mathcal{F}^{0} \rightarrow R^{0}\mathcal{F}$ of functors which is an isomorphism on injective $R$-modules, then there is a unique isomorphism $\tau: (\mathcal{F}^{i})_{i \geq 0} \rightarrow (R^{i}\mathcal{F})_{i \geq 0}$ of $\delta$-functors such that $\tau^{0}=\eta$.
\end{enumerate}
\end{corollary}

\begin{prf}
(i): We note that $(L_{i}\mathcal{F})(F)=0$ for every $i \geq 1$ and every free $R$-module $F$. Now the result follows from Theorem \ref{3.2.4} (i).

(ii): We note that $(R^{i}\mathcal{F})(I)=0$ for every $i \geq 1$ and every injective $R$-module $I$. Now the result follows from Theorem \ref{3.2.4} (ii).
\end{prf}

We next recall the Koszul complex and the Koszul homology briefly. The Koszul complex $K^{R}(a)$ on an element $a \in R$ is the $R$-complex
$$K^{R}(a):=\Cone(R \xrightarrow{a} R),$$
and the Koszul complex $K^{R}(\underline{a})$ on a sequence of elements $\underline{a} = a_{1},\ldots,a_{n} \in R$ is the $R$-complex
$$K^{R}(\underline{a}):= K^{R}(a_{1}) \otimes_{R} \cdots \otimes_{R} K^{R}(a_{n}).$$
It is easy to see that $K^{R}(\underline{a})$ is a complex of finitely generated free $R$-modules concentrated in degrees $n,\ldots,0$.
Given any $R$-module $M$, there is an isomorphism of $R$-complexes
$$K^{R}(\underline{a})\otimes_{R}M \cong \Sigma^{n} \Hom_{R}\left(K^{R}(\underline{a}),M\right),$$
which is sometimes referred to as the self-duality property of the Koszul complex. Accordingly, we feel free to define the Koszul homology of the sequence $\underline{a}$ with coefficients in $M$, by setting
$$H_{i}(\underline{a};M):= H_{i}\left(K^{R}(\underline{a})\otimes_{R}M\right) \cong H_{i-n}\left(\Hom_{R}\left(K^{R}(\underline{a}),M\right)\right)$$
for every $i\geq 0$.

One can form both direct and inverse systems of Koszul complexes and Koszul homologies as explicated in the next remark.

\begin{remark} \label{3.2.6}
We have:
\begin{enumerate}
\item[(i)] Given an element $a \in R$, we define a morphism $\varphi^{k,l}_{a}:K^{R}(a^{k}) \rightarrow K^{R}(a^{l})$ of $R$-complexes for every $k \leq l$ as follows:
\[\begin{tikzpicture}[every node/.style={midway},]
  \matrix[column sep={3em}, row sep={3em}]
  {\node(1) {$0$}; & \node(2) {$R$}; & \node(3) {$R$}; & \node(4) {$0$};\\
  \node(5) {$0$}; & \node(6) {$R$}; & \node(7) {$R$}; & \node(8) {$0$};\\};
  \draw[decoration={markings,mark=at position 1 with {\arrow[scale=1.5]{>}}},postaction={decorate},shorten >=0.5pt] (1) -- (2) node[anchor=east] {};
  \draw[decoration={markings,mark=at position 1 with {\arrow[scale=1.5]{>}}},postaction={decorate},shorten >=0.5pt] (2) -- (3) node[anchor=south] {$a^{k}$};
  \draw[decoration={markings,mark=at position 1 with {\arrow[scale=1.5]{>}}},postaction={decorate},shorten >=0.5pt] (3) -- (4) node[anchor=east] {};
  \draw[decoration={markings,mark=at position 1 with {\arrow[scale=1.5]{>}}},postaction={decorate},shorten >=0.5pt] (5) -- (6) node[anchor=south] {};
  \draw[decoration={markings,mark=at position 1 with {\arrow[scale=1.5]{>}}},postaction={decorate},shorten >=0.5pt] (6) -- (7) node[anchor=south] {$a^{l}$};
  \draw[decoration={markings,mark=at position 1 with {\arrow[scale=1.5]{>}}},postaction={decorate},shorten >=0.5pt] (7) -- (8) node[anchor=south] {};
  \draw[double distance=1.5pt] (2) -- (6) node[anchor=west] {};
  \draw[decoration={markings,mark=at position 1 with {\arrow[scale=1.5]{>}}},postaction={decorate},shorten >=0.5pt] (3) -- (7) node[anchor=west] {$a^{l-k}$};
\end{tikzpicture}\]
It is easily seen that $\left\{K^{R}(a^{k}),\varphi^{k,l}_{a}\right\}_{k \in \mathbb{N}}$ is a direct system of $R$-complexes. Given elements $\underline{a}=a_{1},...,a_{n} \in R$, we let $\underline{a}^{k}=a_{1}^{k},...,a_{n}^{k}$ for every $k \geq 1$. Now
$$K^{R}(\underline{a}^{k})=K^{R}(a^{k}_{1}) \otimes_{R} \cdots \otimes_{R} K^{R}(a^{k}_{n}),$$
and we let
$$\varphi^{k,l}:=\varphi^{k,l}_{a_{1}}\otimes_{R}\cdots\otimes_{R}\varphi^{k,l}_{a_{n}}.$$
It follows that $\left\{K^{R}(\underline{a}^{k}),\varphi^{k,l}\right\}_{k \in \mathbb{N}}$ is a direct system of $R$-complexes. \index{direct system of Koszul complexes} It is also clear that $\left\{H_{i}\left(\underline{a}^{k};M\right),H_{i}\left(\varphi^{k,l}\otimes_{R}M\right)\right\}_{k \in \mathbb{N}}$ is a direct system of $R$-modules for every $i \in \mathbb{Z}$. \index{direct system of Koszul homologies}

\item[(ii)] Given an element $a \in R$, we define a morphism $\psi^{k,l}_{a}:K^{R}(a^{k}) \rightarrow K^{R}(a^{l})$ of $R$-complexes for every $k \geq l$ as follows:
\[\begin{tikzpicture}[every node/.style={midway},]
  \matrix[column sep={3em}, row sep={3em}]
  {\node(1) {$0$}; & \node(2) {$R$}; & \node(3) {$R$}; & \node(4) {$0$};\\
  \node(5) {$0$}; & \node(6) {$R$}; & \node(7) {$R$}; & \node(8) {$0$};\\};
  \draw[decoration={markings,mark=at position 1 with {\arrow[scale=1.5]{>}}},postaction={decorate},shorten >=0.5pt] (1) -- (2) node[anchor=east] {};
  \draw[decoration={markings,mark=at position 1 with {\arrow[scale=1.5]{>}}},postaction={decorate},shorten >=0.5pt] (2) -- (3) node[anchor=south] {$a^{k}$};
  \draw[decoration={markings,mark=at position 1 with {\arrow[scale=1.5]{>}}},postaction={decorate},shorten >=0.5pt] (3) -- (4) node[anchor=east] {};
  \draw[decoration={markings,mark=at position 1 with {\arrow[scale=1.5]{>}}},postaction={decorate},shorten >=0.5pt] (5) -- (6) node[anchor=south] {};
  \draw[decoration={markings,mark=at position 1 with {\arrow[scale=1.5]{>}}},postaction={decorate},shorten >=0.5pt] (6) -- (7) node[anchor=south] {$a^{l}$};
  \draw[decoration={markings,mark=at position 1 with {\arrow[scale=1.5]{>}}},postaction={decorate},shorten >=0.5pt] (7) -- (8) node[anchor=south] {};
  \draw[double distance=1.5pt] (3) -- (7) node[anchor=west] {};
  \draw[decoration={markings,mark=at position 1 with {\arrow[scale=1.5]{>}}},postaction={decorate},shorten >=0.5pt] (2) -- (6) node[anchor=west] {$a^{k-l}$};
\end{tikzpicture}\]
It is easily seen that $\left\{K^{R}(a^{k}),\varphi^{k,l}_{a}\right\}_{k \in \mathbb{N}}$ is an inverse system of $R$-complexes. Given elements $\underline{a}=a_{1},...,a_{n} \in R$, we let $\underline{a}^{k}=a_{1}^{k},...,a_{n}^{k}$ for every $k \geq 1$. Now
$$K^{R}(\underline{a}^{k})=K^{R}(a^{k}_{1}) \otimes_{R} \cdots \otimes_{R} K^{R}(a^{k}_{n}),$$
and we let
$$\psi^{k,l}:=\psi^{k,l}_{a_{1}}\otimes_{R}\cdots\otimes_{R}\psi^{k,l}_{a_{n}}.$$
It follows that $\left\{K^{R}(\underline{a}^{k}),\psi^{k,l}\right\}_{k \in \mathbb{N}}$ is an inverse system of $R$-complexes. \index{inverse system of Koszul complexes} It is also clear that $\left\{H_{i}\left(\underline{a}^{k};M\right),H_{i}\left(\psi^{k,l}\otimes_{R}M\right)\right\}_{k \in \mathbb{N}}$ is an inverse system of $R$-modules for every $i \in \mathbb{Z}$. \index{inverse system of Koszul homologies}
\end{enumerate}
\end{remark}

Recall that an inverse system $\left\{M_{\alpha},\varphi_{\alpha,\beta}\right\}_{\alpha \in \mathbb{N}}$ of $R$-modules is said to satisfy the trivial Mittag-Leffler condition \index{trivial Mittag-Leffler condition} if for every $\beta \in \mathbb{N}$, there is an $\alpha \geq \beta$ such that $\varphi_{\alpha \beta}=0$. Besides, the inverse system $\left\{M_{\alpha},\varphi_{\alpha,\beta}\right\}_{\alpha \in \mathbb{N}}$ of $R$-modules is said to satisfy the Mittag-Leffler condition \index{Mittag-Leffler condition} if for every $\beta \in \mathbb{N}$, there is an $\alpha_{0} \geq \beta$ such that $\im \varphi_{\alpha \beta}=\im \varphi_{\alpha_{0} \beta}$ for every $\alpha \geq \alpha_{0} \geq \beta$. It is straightforward to verify that the trivial Mittag-Leffler condition implies the Mittag-Leffler condition.

The following lemma reveals a significant feature of Koszul homology and lies at the heart of the proof of Greenlees-May Duality. The idea of the proof is taken from \cite{Sc}.

\begin{lemma} \label{3.2.7}
Let $\underline{a}=a_{1},...,a_{n} \in R$, and $\underline{a}^{k}=a_{1}^{k},...,a_{n}^{k}$ for every $k \geq 1$. Then the inverse system $\left\{H_{i}\left(\underline{a}^{k};R\right)\right\}_{k \in \mathbb{N}}$ satisfies the trivial Mittag-Leffler condition for every $i \geq 1$.
\end{lemma}

\begin{prf}
Let $a\in R$ and $M$ a finitely generated $R$-module. The transition maps of the inverse system $\left\{K^{R}(a^{k})\otimes_{R}M\right\}_{k \in \mathbb{N}}$ can be identified with the following morphisms of $R$-complexes for every $k \geq l$:
\[\begin{tikzpicture}[every node/.style={midway},]
  \matrix[column sep={3em}, row sep={3em}]
  {\node(1) {$0$}; & \node(2) {$M$}; & \node(3) {$M$}; & \node(4) {$0$};\\
  \node(5) {$0$}; & \node(6) {$M$}; & \node(7) {$M$}; & \node(8) {$0$};\\};
  \draw[decoration={markings,mark=at position 1 with {\arrow[scale=1.5]{>}}},postaction={decorate},shorten >=0.5pt] (1) -- (2) node[anchor=east] {};
  \draw[decoration={markings,mark=at position 1 with {\arrow[scale=1.5]{>}}},postaction={decorate},shorten >=0.5pt] (2) -- (3) node[anchor=south] {$a^{k}$};
  \draw[decoration={markings,mark=at position 1 with {\arrow[scale=1.5]{>}}},postaction={decorate},shorten >=0.5pt] (3) -- (4) node[anchor=east] {};
  \draw[decoration={markings,mark=at position 1 with {\arrow[scale=1.5]{>}}},postaction={decorate},shorten >=0.5pt] (5) -- (6) node[anchor=south] {};
  \draw[decoration={markings,mark=at position 1 with {\arrow[scale=1.5]{>}}},postaction={decorate},shorten >=0.5pt] (6) -- (7) node[anchor=south] {$a^{l}$};
  \draw[decoration={markings,mark=at position 1 with {\arrow[scale=1.5]{>}}},postaction={decorate},shorten >=0.5pt] (7) -- (8) node[anchor=south] {};
  \draw[double distance=1.5pt] (3) -- (7) node[anchor=west] {};
  \draw[decoration={markings,mark=at position 1 with {\arrow[scale=1.5]{>}}},postaction={decorate},shorten >=0.5pt] (2) -- (6) node[anchor=west] {$a^{k-l}$};
\end{tikzpicture}\]
Since $H_{1}\left(a^{k};M\right)=\left(0:_{M}a^{k}\right)$, the transition maps of the inverse system $\left\{H_{1}\left(a^{k};M\right)\right\}_{k \in \mathbb{N}}$ can be identified with the $R$-homomorphisms
$$\left(0:_{M}a^{k}\right) \xrightarrow{a^{k-l}} \left(0:_{M}a^{l}\right)$$
for every $k \geq l$. Fix $l \in \mathbb{N}$. Since $R$ is noetherian and $M$ is finitely generated, the ascending chain
$$\left(0:_{M}a\right) \subseteq \left(0:_{M}a^{2}\right) \subseteq \cdots$$
of submodules of $M$ stabilizes, i.e. there is an integer $t \geq 1$ such that
$$\left(0:_{M}a^{t}\right) = \left(0:_{M}a^{t+1}\right) = \cdots.$$
Set $k:=t+l$. Then the transition map $\left(0:_{M}a^{k}\right) \xrightarrow{a^{k-l}} \left(0:_{M}a^{l}\right)$ is zero. Indeed, if $x \in \left(0:_{M}a^{k}\right)$, then since
$$\left(0:_{M}a^{k}\right) = \left(0:_{M}a^{t+l}\right) = \left(0:_{M}a^{t}\right),$$
we have $x \in \left(0:_{M}a^{t}\right)$, so $a^{k-l}x=a^{t}x=0$. This shows that the inverse system $\left\{H_{1}\left(a^{k};M\right)\right\}_{k \in \mathbb{N}}$ satisfies the trivial Mittag-Leffler condition. But $H_{i}\left(a^{k};M\right)=0$ for every $i \geq 2$, so the inverse system $\left\{H_{i}\left(a^{k};M\right)\right\}_{k \in \mathbb{N}}$ satisfies the trivial Mittag-Leffler condition for every $i \geq 1$.

Now we argue by induction on $n$. If $n=1$, then the inverse system $\left\{H_{i}\left(a_{1}^{k};R\right)\right\}_{k \in \mathbb{N}}$ satisfies the trivial Mittag-Leffler condition for every $i \geq 1$ by the discussion above. Now assume that $n \geq 2$, and make the obvious induction hypothesis. There is an exact sequence of inverse systems
\begin{equation} \label{eq:3.2.7.1}
\left\{H_{i}\left(a_{1}^{k},...,a_{n-1}^{k};R\right)\right\}_{k \in \mathbb{N}} \rightarrow \left\{H_{0}\left(a_{n}^{k};H_{i}\left(a_{1}^{k},...,a_{n-1}^{k};R\right)\right)\right\}_{k \in \mathbb{N}} \rightarrow 0
\end{equation}
of $R$-modules for every $i \geq 0$. By the induction hypothesis, the inverse system $\left\{H_{i}\left(a_{1}^{k},...,a_{n-1}^{k};R\right)\right\}_{k \in \mathbb{N}}$ satisfies the trivial Mittag-Leffler condition for every $i \geq 1$, so the exact sequence \eqref{eq:3.2.7.1} shows that the inverse system $\left\{H_{0}\left(a_{n}^{k};H_{i}\left(a_{1}^{k},...,a_{n-1}^{k};R\right)\right)\right\}_{k \in \mathbb{N}}$ satisfies the Mittag-Leffler condition for every $i \geq 1$. On the other hand, there is a short exact sequence of inverse systems
\begin{equation} \label{eq:3.2.7.2}
0 \rightarrow \left\{H_{0}\left(a_{n}^{k};H_{i}\left(a_{1}^{k},...,a_{n-1}^{k};R\right)\right)\right\}_{k \in \mathbb{N}} \rightarrow \left\{H_{i}\left(a_{1}^{k},...,a_{n}^{k};R\right)\right\}_{k \in \mathbb{N}} \rightarrow $$$$ \left\{H_{1}\left(a_{n}^{k};H_{i-1}\left(a_{1}^{k},...,a_{n-1}^{k};R\right)\right)\right\}_{k \in \mathbb{N}} \rightarrow 0
\end{equation}
of $R$-modules for every $i \geq 0$.
Since $H_{i-1}\left(a_{1}^{k},...,a_{n-1}^{k};R\right)$ is a finitely generated $R$-module for every $i \geq 1$, the discussion above shows that
$\left\{H_{1}\left(a_{n}^{k};H_{i-1}\left(a_{1}^{k},...,a_{n-1}^{k};R\right)\right)\right\}_{k \in \mathbb{N}}$ satisfies the Mittag-Leffler condition for every $i \geq 1$. Therefore, the short exact sequence \eqref{eq:3.2.7.2} shows that the inverse system $\left\{H_{i}\left(a_{1}^{k},...,a_{n}^{k};R\right)\right\}_{k \in \mathbb{N}}$ satisfies the trivial Mittag-Leffler condition for every $i \geq 1$.
\end{prf}

The category $\mathcal{C}(R)$ of $R$-complexes enjoys direct limits and inverse limits. However, the derived category $\mathcal{D}(R)$ does not support the notions of direct limits and inverse limits. But this situation is remedied by the existence of homotopy direct limits and homotopy inverse limits as defined in triangulated categories with countable products and coproducts.

\begin{remark} \label{3.2.8}
Let $\left\{X^{\alpha},\varphi^{\alpha \beta}\right\}_{\alpha \in \mathbb{N}}$ be a direct system of $R$-complexes, and $\left\{Y^{\alpha},\psi^{\alpha \beta}\right\}_{\alpha \in \mathbb{N}}$ an inverse system of $R$-complexes. Then we have:
\begin{enumerate}
\item[(i)] The direct limit \index{direct limit of complexes} of the direct system $\left\{X^{\alpha},\varphi^{\alpha \beta}\right\}_{\alpha \in \mathbb{N}}$ is an $R$-complex $\varinjlim X^{\alpha}$ given by $\left(\varinjlim X^{\alpha}\right)_{i}=\varinjlim X^{\alpha}_{i}$ and $\partial^{\varinjlim X^{\alpha}}_{i}=\varinjlim \partial^{X^{\alpha}}_{i}$ for every $i \in \mathbb{Z}$. Indeed, it is easy to see that $\varinjlim X^{\alpha}$ satisfies the universal property of direct limits in a category.
\item[(ii)] The homotopy direct limit \index{homotopy direct limit of complexes} of the direct system $\left\{X^{\alpha},\varphi^{\alpha \beta}\right\}_{\alpha \in \mathbb{N}}$ is given by $\hocolim X^{\alpha} = \Cone(\vartheta)$, where the morphism $\vartheta: \bigoplus_{\alpha=1}^{\infty}X^{\alpha} \rightarrow \bigoplus_{\alpha=1}^{\infty}X^{\alpha}$ is given by $\vartheta_{i}\left((x_{i}^{\alpha})\right)= \iota_{i}^{\alpha}(x_{i}^{\alpha})-\iota_{i}^{\alpha+1}\left(\varphi_{i}^{\alpha,\alpha+1}(x_{i}^{\alpha})\right)$ for every $i \in \mathbb{Z}$. Indeed, it is easy to see that the morphism $\vartheta$ fits into a distinguished triangle
    $$\bigoplus_{\alpha=1}^{\infty}X^{\alpha} \rightarrow \bigoplus_{\alpha=1}^{\infty}X^{\alpha} \rightarrow \hocolim X^{\alpha} \rightarrow.$$
\item[(iii)] The inverse limit \index{inverse limit of complexes} of the inverse system $\left\{Y^{\alpha},\psi^{\alpha \beta}\right\}_{\alpha \in \mathbb{N}}$ is an $R$-complex $\varprojlim Y^{\alpha}$ given by $\left(\varprojlim Y^{\alpha}\right)_{i}=\varprojlim Y^{\alpha}_{i}$ and $\partial^{\varprojlim Y^{\alpha}}_{i}=\varprojlim \partial^{Y^{\alpha}}_{i}$ for every $i \in \mathbb{Z}$. Indeed, it is easy to see that $\varprojlim X^{\alpha}$ satisfies the universal property of inverse limits in a category.
\item[(iv)] The homotopy inverse limit \index{homotopy inverse limit of complexes} of the inverse system $\left\{Y^{\alpha},\psi^{\alpha \beta}\right\}_{\alpha \in \mathbb{N}}$ is given by $\holim Y^{\alpha} = \Sigma^{-1}\Cone(\varpi)$, where the morphism $\varpi: \prod_{\alpha=1}^{\infty}Y^{\alpha} \rightarrow \prod_{\alpha=1}^{\infty}Y^{\alpha}$ is given by $\varpi_{i}\left((y_{i}^{\alpha})\right)= \left(y_{i}^{\alpha}-\psi_{i}^{\alpha+1,\alpha}(y_{i}^{\alpha+1})\right)$ for every $i \in \mathbb{Z}$. Indeed, it is easy to see that the morphism $\varpi$ fits into a distinguished triangle
    $$\holim Y^{\alpha} \rightarrow \prod_{\alpha=1}^{\infty}Y^{\alpha} \rightarrow \prod_{\alpha=1}^{\infty}Y^{\alpha} \rightarrow.$$
\end{enumerate}
\end{remark}

The Mittag-Leffler condition forces many limits to be zero.

\begin{lemma} \label{3.2.9}
Let $\left\{M_{\alpha},\varphi_{\alpha \beta}\right\}_{\alpha \in \mathbb{N}}$ be an inverse system of $R$-modules that satisfies the trivial Mittag-Leffler condition, and $\mathcal{F}:\mathcal{M}(R) \rightarrow \mathcal{M}(R)$ an additive contravariant functor. Then the following assertions hold:
\begin{enumerate}
\item[(i)] $\varprojlim M_{\alpha} = 0 = \varprojlim^{1} M_{\alpha}$.
\item[(ii)] $\varinjlim \mathcal{F}(M_{\alpha})=0$.
\end{enumerate}
\end{lemma}

\begin{prf}
(i): Let $\varpi:\prod_{\alpha \in \mathbb{N}}M_{\alpha} \rightarrow \prod_{\alpha \in \mathbb{N}}M_{\alpha}$ be an $R$-homomorphism given by $\varpi\left((x_{\alpha})\right)=\left(x_{\alpha}-\varphi_{\alpha+1,\alpha}(x_{\alpha+1})\right)$. We show that $\varpi$ is an isomorphism. Let $(x_{\alpha}) \in \prod_{\alpha \in \mathbb{N}}M_{\alpha}$ be such that $x_{\alpha}=\varphi_{\alpha+1,\alpha}(x_{\alpha+1})$ for every $\alpha \in \mathbb{N}$. Fix $\alpha \in \mathbb{N}$, and by the trivial Mittag-Leffler condition choose $\gamma \geq \alpha$ such that $\varphi_{\gamma \alpha}=0$. Then we have
\begin{equation*}
\begin{split}
 x_{\alpha} & = \varphi_{\alpha+1,\alpha}(x_{\alpha+1}) \\
 & = \varphi_{\alpha+1,\alpha}\left(\varphi_{\alpha+2,\alpha+1}\left(...\left(\varphi_{\gamma,\gamma-1} (x_{\gamma})\right)\right)\right)\\
 & = \varphi_{\gamma \alpha}(x_{\gamma})\\
 & = 0.\\
\end{split}
\end{equation*}
Hence $(x_{\alpha})=0$, and thus $\varpi$ is injective. Now let $(y_{\alpha}) \in \prod_{\alpha \in \mathbb{N}}M_{\alpha}$. For any $\beta \in \mathbb{N}$, we set $x_{\beta}:= \sum_{\alpha=\beta}^{\infty}\varphi_{\alpha \beta}(y_{\alpha})$ which is a finite sum by the trivial Mittag-Leffler condition. Then we have
\begin{equation*}
\begin{split}
 \varphi_{\beta+1,\beta}(x_{\beta+1}) & = \varphi_{\beta+1,\beta}\left(\sum_{\alpha=\beta+1}^{\infty}\varphi_{\alpha, \beta+1}(y_{\alpha})\right) \\
 & = \sum_{\alpha=\beta+1}^{\infty} \varphi_{\alpha \beta}(y_{\alpha})\\
 & = \sum_{\alpha=\beta}^{\infty} \varphi_{\alpha \beta}(y_{\alpha})-\varphi_{\beta \beta}(y_{\beta})\\
 & = x_{\beta}-y_{\beta}.\\
\end{split}
\end{equation*}
Therefore, we have
$$\varpi\left((x_{\alpha})\right)=\left(x_{\alpha}-\varphi_{\alpha+1,\alpha}(x_{\alpha+1})\right)=(y_{\alpha}),$$
so $\varpi$ is surjective. It follows that $\varpi$ is an isomorphism. Therefore, $\varprojlim M_{\alpha} \cong \ker \varpi=0$ and $\varprojlim^{1} M_{\alpha} \cong \coker \varpi=0$.

(ii): First we note that $\left\{\mathcal{F}(M_{\alpha}),\psi_{\beta\alpha}:=\mathcal{F}(\varphi_{\alpha\beta})\right\}_{\alpha \in \mathbb{N}}$ is a direct system of $R$-modules. Let $\psi_{\alpha}:\mathcal{F}(M_{\alpha}) \rightarrow \varinjlim \mathcal{F}(M_{\alpha})$ be the natural injection of direct limit for every $\alpha \in \mathbb{N}$. We know that an arbitrary element of $\varinjlim \mathcal{F}(M_{\alpha})$ is of the form $\psi_{t}(y)$ for some $t \in \mathbb{N}$ and some $y \in \mathcal{F}(M_{t})$. By the trivial Mittag-Leffler condition, there is an integer $s \geq t$ such that $\varphi_{st}=0$, so that $\psi_{ts}=\mathcal{F}(\varphi_{st})=0$. Then $\psi_{t}(y)=\psi_{s}\left(\psi_{ts}(y)\right)=0$. Hence $\varinjlim \mathcal{F}(M_{\alpha})=0$.
\end{prf}

The next proposition collects some information on the homology of limits.

\begin{proposition} \label{3.2.10}
Let $\left\{X^{\alpha},\varphi^{\alpha \beta}\right\}_{\alpha \in \mathbb{N}}$ be a direct system of $R$-complexes, and $\left\{Y^{\alpha},\psi^{\alpha \beta}\right\}_{\alpha \in \mathbb{N}}$ an inverse system of $R$-complexes. Then the following assertions hold for every $i \in \mathbb{Z}$:
\begin{enumerate}
\item[(i)] There is a natural isomorphism $H_{i}\left(\varinjlim X^{\alpha}\right) \cong \varinjlim H_{i}(X^{\alpha})$.
\item[(ii)] There is a natural isomorphism $H_{i}\left(\hocolim X^{\alpha}\right) \cong \varinjlim H_{i}(X^{\alpha})$.
\item[(iii)] If the inverse system $\left\{Y^{\alpha}_{i},\psi^{\alpha \beta}_{i}\right\}_{\alpha \in \mathbb{N}}$ of $R$-modules satisfies the Mittag-Leffler condition for every $i \in \mathbb{Z}$, then there is a short exact sequence
\begin{center}
$0 \rightarrow \varprojlim^{1} H_{i+1}(Y^{\alpha}) \rightarrow H_{i}\left(\varprojlim Y^{\alpha}\right) \rightarrow \varprojlim H_{i}(Y^{\alpha}) \rightarrow 0$
\end{center}
of $R$-modules.
\item[(iv)] There is a short exact sequence
\begin{center}
$0 \rightarrow \varprojlim^{1} H_{i+1}(Y^{\alpha}) \rightarrow H_{i}\left(\holim Y^{\alpha}\right) \rightarrow \varprojlim H_{i}(Y^{\alpha}) \rightarrow 0$
\end{center}
of $R$-modules.
\end{enumerate}
\end{proposition}

\begin{prf}
(i): See \cite[Theorem 4.2.4]{Se}.

(ii): See the paragraph before \cite[Lemma 0.1]{GM}.

(iii): See \cite[Theorem 3.5.8]{We}.

(iv): See the paragraph after \cite[Lemma 0.1]{GM}.
\end{prf}

Now we are ready to present the following definitions.

\begin{definition} \label{3.2.11}
Let $\underline{a}=a_{1},...,a_{n} \in R$. Then:
\begin{enumerate}
\item[(i)] Define the \textit{\v{C}ech complex} \index{\v{C}ech complex} on the elements $\underline{a}$ to be $\check{C}(\underline{a}):= \varinjlim \Sigma^{-n}K^{R}(\underline{a}^{k})$.
\item[(ii)] Define the \textit{stable \v{C}ech complex} \index{stable \v{C}ech complex} on the elements $\underline{a}$ to be $\check{C}_{\infty}(\underline{a}):= \hocolim \Sigma^{-n}K^{R}(\underline{a}^{k})$.
\end{enumerate}
\end{definition}

We note that $\check{C}(\underline{a})$ is a bounded $R$-complex of flat modules concentrated in degrees $0,-1,...,-n$, and $\check{C}_{\infty}(\underline{a})$ is a bounded $R$-complex of free modules concentrated in degrees $1,0,...,-n$. Moreover, it can be shown that there is a quasi-isomorphism $\check{C}_{\infty}(\underline{a}) \xrightarrow{\simeq} \check{C}(\underline{a})$, which in turn implies that $\check{C}_{\infty}(\underline{a}) \simeq \check{C}(\underline{a})$ in $\mathcal{D}(R)$. Therefore, $\check{C}_{\infty}(\underline{a})$ is a semi-projective approximation of the semi-flat $R$-complex $\check{C}(\underline{a})$.

The next proposition investigates the relation between local cohomology and local homology with \v{C}ech complex and stable \v{C}ech complex, and provides the first essential step towards the Greenlees-May Duality.

\begin{proposition} \label{3.2.12}
Let $\mathfrak{a}= (a_{1},...,a_{n})$ be an ideal of $R$, $\underline{a}=a_{1},...,a_{n}$, and $M$ an $R$-module. Then there are natural isomorphisms for every $i \geq 0$:
\begin{enumerate}
\item[(i)] $H^{i}_{\mathfrak{a}}(M) \cong H_{-i}\left(\check{C}(\underline{a})\otimes_{R}M\right) \cong H_{-i}\left(\check{C}_{\infty}(\underline{a})\otimes_{R}M\right)$.
\item[(ii)] $H^{\mathfrak{a}}_{i}(M) \cong H_{i}\left(\Hom_{R}\left(\check{C}_{\infty}(\underline{a}),M\right)\right)$.
\end{enumerate}
\end{proposition}

\begin{prf}
(i): Let $\mathcal{F}^{i}=H_{-i}\left(\check{C}(\underline{a})\otimes_{R}-\right):\mathcal{M}(R) \rightarrow \mathcal{M}(R)$ for every $i \geq 0$. Given a short exact sequence
$$0 \rightarrow M' \rightarrow M \rightarrow M'' \rightarrow 0$$
of $R$-modules, since $\check{C}(\underline{a})$ is an $R$-complex of flat modules, the functor $\check{C}(\underline{a})\otimes_{R}-:\mathcal{C}(R) \rightarrow \mathcal{C}(R)$ is exact, whence we get a short exact sequence
$$0 \rightarrow \check{C}(\underline{a})\otimes_{R}M' \rightarrow \check{C}(\underline{a})\otimes_{R}M \rightarrow \check{C}(\underline{a})\otimes_{R}M'' \rightarrow 0$$
of $R$-complexes, which in turn yields a long exact homology sequence in a functorial way. This shows that $\left(\mathcal{F}^{i}:\mathcal{M}(R) \rightarrow \mathcal{M}(R)\right)_{i \geq 0}$ is a cohomological covariant $\delta$-functor. Moreover, using Proposition \ref{3.2.10} (i), we have
\begin{equation} \label{eq:3.2.12.1}
\begin{split}
 \mathcal{F}^{i} & = H_{-i}\left(\check{C}(\underline{a})\otimes_{R}-\right) \\
 & = H_{-i}\left(\left(\varinjlim \Sigma^{-n}K^{R}(\underline{a}^{k})\right)\otimes_{R}-\right) \\
 & \cong \varinjlim H_{n-i}\left(K^{R}(\underline{a}^{k})\otimes_{R}-\right) \\
 & \cong \varinjlim H_{n-i}\left(\underline{a}^{k};-\right)\\
\end{split}
\end{equation}
for every $i \geq 0$.

Let $I$ be an injective $R$-module. Then by the display \eqref{eq:3.2.12.1}, we have
\begin{equation} \label{eq:3.2.12.2}
\begin{split}
 \mathcal{F}^{i}(I) & = \varinjlim H_{n-i}\left(\underline{a}^{k};I\right)\\
 & \cong \varinjlim H_{-i}\left(\Hom_{R}\left(K^{R}(\underline{a}^{k}),I\right)\right) \\
 & \cong \varinjlim \Hom_{R}\left(H_{i}\left(K^{R}(\underline{a}^{k})\right),I\right) \\
 & \cong \varinjlim \Hom_{R}\left(H_{i}\left(\underline{a}^{k};R\right),I\right).\\
\end{split}
\end{equation}
By Lemma \ref{3.2.7}, the inverse system $\left\{H_{i}\left(\underline{a}^{k};R\right)\right\}_{k \in \mathbb{N}}$ satisfies the trivial Mittag-Leffler condition for every $i \geq 1$. Now Lemma \ref{3.2.9} (ii) implies that $\varinjlim \Hom_{R}\left(H_{i}\left(\underline{a}^{k};R\right),I\right)=0$, thereby the display \eqref{eq:3.2.12.2} shows that $\mathcal{F}^{i}(I)=0$ for every $i \geq 1$.

Let $M$ be an $R$-module. Then by the display \eqref{eq:3.2.12.1}, we have the natural isomorphisms
\begin{equation*}
\begin{split}
 \mathcal{F}^{0}(M) & \cong \varinjlim H_{n}\left(\underline{a}^{k};M\right) \\
 & \cong \varinjlim \left(0:_{M}(\underline{a}^{k})\right) \\
 & \cong \varinjlim \Hom_{R}\left(R/(\underline{a}^{k}),M\right) \\
 & \cong \varinjlim \Hom_{R}\left(R/\mathfrak{a}^{k},M\right)\\
 & \cong \Gamma_{\mathfrak{a}}(M)\\
 & \cong H^{0}_{\mathfrak{a}}(M).\\
\end{split}
\end{equation*}
It follows from Corollary \ref{3.2.5} (ii) that $H^{i}_{\mathfrak{a}}(-) \cong \mathcal{F}^{i}$ for every $i \geq 0$.

For the second isomorphism, using the display \eqref{eq:3.2.12.1} and Proposition \ref{3.2.10} (ii), we have the natural isomorphisms
\begin{equation*}
\begin{split}
 H^{i}_{\mathfrak{a}}(M) & \cong \mathcal{F}^{i}(M) \\
 & \cong \varinjlim H_{n-i}\left(\underline{a}^{k};M\right) \\
 & \cong \varinjlim H_{n-i}\left(K^{R}(\underline{a}^{k})\otimes_{R}M\right) \\
 & \cong H_{n-i}\left(\hocolim\left(K^{R}(\underline{a}^{k})\otimes_{R}M\right)\right) \\
 & \cong H_{-i}\left(\left(\hocolim \Sigma^{-n}K^{R}(\underline{a})\right)\otimes_{R}M\right) \\
 & \cong H_{-i}\left(\check{C}_{\infty}(\underline{a})\otimes_{R}M\right) \\
\end{split}
\end{equation*}
for every $i \geq 0$.

(ii): Let $\mathcal{F}_{i}=H_{i}\left(\Hom_{R}\left(\check{C}_{\infty}(\underline{a}),-\right)\right):\mathcal{M}(R) \rightarrow \mathcal{M}(R)$ for every $i \geq 0$. Given a short exact sequence
\[0 \rightarrow M' \rightarrow M \rightarrow M'' \rightarrow 0\]
$R$-modules, since $\check{C}_{\infty}(\underline{a})$ is an $R$-complex of free modules, the functor $\Hom_{R}\left(\check{C}_{\infty}(\underline{a}),-\right):\mathcal{C}(R) \rightarrow \mathcal{C}(R)$ is exact, whence we get a short exact sequence
$$0 \rightarrow \Hom_{R}\left(\check{C}_{\infty}(\underline{a}),M'\right) \rightarrow \Hom_{R}\left(\check{C}_{\infty}(\underline{a}),M\right) \rightarrow \Hom_{R}\left(\check{C}_{\infty}(\underline{a}),M''\right) \rightarrow 0$$
of $R$-complexes, which in turn yields a long exact homology sequence in a functorial way. It follows that $\left(\mathcal{F}_{i}:\mathcal{M}(R) \rightarrow \mathcal{M}(R)\right)_{i \geq 0}$ is a homological covariant $\delta$-functor. Moreover, using the self-duality property of Koszul complex, we have
\begin{equation} \label{eq:3.2.12.3}
\begin{split}
 \mathcal{F}_{i} & = H_{i}\left(\Hom_{R}\left(\check{C}_{\infty}(\underline{a}),-\right)\right) \\
 & = H_{i}\left(\Hom_{R}\left(\hocolim \Sigma^{-n}K^{R}(\underline{a}^{k}),-\right)\right) \\
 & \cong H_{i}\left(\holim \Sigma^{n}\Hom_{R}\left(K^{R}(\underline{a}^{k}),-\right)\right)\\
 & \cong H_{i}\left(\holim \left(K^{R}(\underline{a}^{k})\otimes_{R}-\right)\right)\\
\end{split}
\end{equation}
for every $i \geq 0$.

Let $M$ be an $R$-module. By Proposition \ref{3.2.10} (iv), we get a short exact sequence
\begin{center}
$0 \rightarrow \varprojlim ^{1} H_{i+1}\left(K^{R}(\underline{a}^{k})\otimes_{R}M\right) \rightarrow H_{i}\left(\holim\left(K^{R}(\underline{a}^{k})\otimes_{R}M\right)\right) \rightarrow \varprojlim H_{i}\left(K^{R}(\underline{a}^{k})\otimes_{R}M\right) \rightarrow 0,$
\end{center}
which implies the short exact sequence
\begin{center}
$0 \rightarrow \varprojlim^{1} H_{i+1}\left(\underline{a}^{k};M\right) \rightarrow \mathcal{F}_{i}(M) \rightarrow \varprojlim H_{i}\left(\underline{a}^{k};M\right) \rightarrow 0$
\end{center}
of $R$-modules for every $i \geq 0$.

Let $F$ be a free $R$-module. If $i \geq 1$, then the inverse system $\left\{H_{i}\left(\underline{a}^{k};R\right)\right\}_{k \in \mathbb{N}}$ satisfies the trivial Mittag-Leffler condition by Lemma \ref{3.2.7}. But $H_{i}\left(\underline{a}^{k};F\right) \cong H_{i}\left(\underline{a}^{k};R\right)\otimes_{R}F$, so it straightforward to see that the inverse system $\left\{H_{i}\left(\underline{a}^{k};F\right)\right\}_{k \in \mathbb{N}}$ satisfies the trivial Mittag-Leffler condition for every $i \geq 1$. Therefore, Lemma \ref{3.2.9} (i) implies that
\begin{center}
$\varprojlim^{1} H_{i}\left(\underline{a}^{k};F\right)=0=\varprojlim H_{i}\left(\underline{a}^{k};F\right)$
\end{center}
for every $i \geq 1$. It follows from the above short exact sequence that $\mathcal{F}_{i}(F)=0$ for every $i \geq 1$.

Upon setting $i=0$, the above short exact sequence yields
\begin{center}
$0 = \varprojlim^{1} H_{1}\left(\underline{a}^{k};F\right) \rightarrow \mathcal{F}_{0}(F) \rightarrow \varprojlim H_{0}\left(\underline{a}^{k};F\right) \rightarrow 0.$
\end{center}
Thus we get the natural isomorphisms
\begin{equation*}
\begin{split}
 \mathcal{F}_{0}(F) & \cong \varprojlim H_{0}\left(\underline{a}^{k};F\right) \\
 & \cong \varprojlim F/(\underline{a}^{k})F \\
 & \cong \varprojlim F/\mathfrak{a}^{k}F \\
 & = \widehat{F}^{\mathfrak{a}}\\
 & \cong H^{\mathfrak{a}}_{0}(F).\\
\end{split}
\end{equation*}
It now follows from Corollary \ref{3.2.5} (i) that $H^{\mathfrak{a}}_{i}(-) \cong \mathcal{F}_{i}$ for every $i \geq 0$.
\end{prf}

\begin{remark} \label{3.2.13}
One should note that $H^{\mathfrak{a}}_{i}(M) \ncong H_{i}\left(\Hom_{R}\left(\check{C}(\underline{a}),M\right)\right)$.
\end{remark}

\section{Complex Prerequisites}

In this section, we commence on developing the requisite tools on complexes which are to be deployed in Section 4. For more information on the material in this section, refer to \cite{AF}, \cite{Ha}, \cite{Fo}, \cite{Li}, and \cite{Sp}.

The derived category $\mathcal{D}(R)$ \index{derived category} is defined as the localization of the homotopy category $\mathcal{K}(R)$ with respect to the multiplicative system of quasi-isomorphisms. Simply put, an object in $\mathcal{D}(R)$ is an $R$-complex $X$ displayed in the standard homological style
$$X= \cdots \rightarrow X_{i+1} \xrightarrow {\partial^{X}_{i+1}} X_{i} \xrightarrow {\partial^{X}_{i}} X_{i-1} \rightarrow \cdots,$$
and a morphism $\varphi:X\rightarrow Y$ in $\mathcal{D}(R)$ is given by the equivalence class of a pair $(f,g)$ of morphisms
$X \xleftarrow{g} U \xrightarrow{f} Y$ in $\mathcal{C}(R)$ with $g$ a quasi-isomorphism, under the equivalence relation that identifies two such pairs $(f,g)$ and $(f^{\prime},g^{\prime})$, whenever there is a diagram in $\mathcal{C}(R)$ as follows which commutes up to homotopy:
\[\begin{tikzpicture}[every node/.style={midway}]
  \matrix[column sep={3em}, row sep={3em}]
  {\node(1) {$$}; & \node(2) {$U$}; & \node(3) {$$};\\
  \node(4) {$X$}; & \node(5) {$V$}; & \node(6) {$Y$};\\
  \node(7) {$$}; & \node(8) {$U^{\prime}$}; & \node(9) {$$};\\};
  \draw[decoration={markings,mark=at position 1 with {\arrow[scale=1.5]{>}}},postaction={decorate},shorten >=0.5pt] (2) -- (4) node[anchor=south] {$g$} node[anchor=west] {$\simeq$};
  \draw[decoration={markings,mark=at position 1 with {\arrow[scale=1.5]{>}}},postaction={decorate},shorten >=0.5pt] (2) -- (6) node[anchor=south] {$f$};
  \draw[decoration={markings,mark=at position 1 with {\arrow[scale=1.5]{>}}},postaction={decorate},shorten >=0.5pt] (8) -- (4) node[anchor=north] {$g^{\prime}$} node[anchor=west] {$\simeq$};
  \draw[decoration={markings,mark=at position 1 with {\arrow[scale=1.5]{>}}},postaction={decorate},shorten >=0.5pt] (8) -- (6) node[anchor=north] {$f^{\prime}$};
  \draw[decoration={markings,mark=at position 1 with {\arrow[scale=1.5]{>}}},postaction={decorate},shorten >=0.5pt] (5) -- (2) node[anchor=east] {$$};
  \draw[decoration={markings,mark=at position 1 with {\arrow[scale=1.5]{>}}},postaction={decorate},shorten >=0.5pt] (5) -- (4) node[anchor=south] {$\simeq$};
  \draw[decoration={markings,mark=at position 1 with {\arrow[scale=1.5]{>}}},postaction={decorate},shorten >=0.5pt] (5) -- (6) node[anchor=east] {$$};
  \draw[decoration={markings,mark=at position 1 with {\arrow[scale=1.5]{>}}},postaction={decorate},shorten >=0.5pt] (5) -- (8) node[anchor=east] {$$};
\end{tikzpicture}\]
The isomorphisms in $\mathcal{D}(R)$ are marked by the symbol $\simeq$.

The derived category $\mathcal{D}(R)$ is triangulated. A distinguished triangle \index{distinguished triangle} in $\mathcal{D}(R)$ is a triangle that is isomorphic to a triangle of the form
$$X \xrightarrow {\mathfrak{L}(f)} Y \xrightarrow{\mathfrak{L}(\varepsilon)} \Cone(f) \xrightarrow{\mathfrak{L}(\varpi)} \Sigma X,$$
for some morphism $f:X \rightarrow Y$ in $\mathcal{C}(R)$ with the mapping cone sequence
$$0 \rightarrow Y \xrightarrow{\varepsilon} \Cone(f) \xrightarrow{\varpi} \Sigma X \rightarrow 0,$$
in which $\mathfrak{L}:\mathcal{C}(R) \rightarrow \mathcal{D}(R)$ is the canonical functor that is defined as $\mathfrak{L}(X)=X$ for every $R$-complex $X$, and $\mathfrak{L}(f)=\varphi$ where $\varphi$ is represented by the morphisms $X \xleftarrow{1^{X}} X \xrightarrow{f} Y$ in $\mathcal{C}(R)$. We note that if $f$ is a quasi-isomorphism in $\mathcal{C}(R)$, then $\mathfrak{L}(f)$ is an isomorphism in $\mathcal{D}(R)$. We sometimes use the shorthand notation
$$X \rightarrow Y \rightarrow Z \rightarrow$$
for a distinguished triangle.

We let $\mathcal{D}_{\sqsubset}(R)$ (res. $\mathcal{D}_{\sqsupset}(R)$) denote the full subcategory of $\mathcal{D}(R)$ consisting of $R$-complexes $X$ with $H_{i}(X)=0$ for $i \gg 0$ (res. $i \ll 0$), and  $D_{\square}(R):=\mathcal{D}_{\sqsubset}(R)\cap \mathcal{D}_{\sqsupset}(R)$. We further let $\mathcal{D}^{f}(R)$ denote the full subcategory of $\mathcal{D}(R)$ consisting of $R$-complexes $X$ with finitely generated homology modules. We also feel free to use any combination of the subscripts and the superscript as in $\mathcal{D}^{f}_{\square}(R)$, with the obvious meaning of the intersection of the two subcategories involved.

We recall the resolutions of complexes.

\begin{definition} \label{3.3.1}
We have:
\begin{enumerate}
\item[(i)] An $R$-complex $P$ of projective modules is said to be \textit{semi-projective} \index{semi-projective complex} if the functor $\Hom_{R}(P,-)$ preserves quasi-isomorphisms. By a \textit{semi-projective resolution} \index{semi-projective resolution} of an $R$-complex $X$, we mean a quasi-isomorphism $P\xrightarrow {\simeq} X$ in which $P$ is a semi-projective $R$-complex.
\item[(ii)] An $R$-complex
$I$ of injective modules is said to be \textit{semi-injective} \index{semi-injective complex} if the functor $\Hom_{R}(-,I)$ preserves quasi-isomorphisms. By a \textit{semi-injective resolution} \index{semi-injective resolution} of an $R$-complex $X$, we mean a quasi-isomorphism $X\xrightarrow {\simeq} I$ in which $I$ is a semi-injective $R$-complex.
\item[(iii)] An $R$-complex
$F$ of flat modules is said to be \textit{semi-flat} \index{semi-flat complex} if the functor $F\otimes_{R}-$ preserves quasi-isomorphisms. By a \textit{semi-flat resolution} \index{semi-flat resolution} of an $R$-complex $X$, we mean a quasi-isomorphism $F\xrightarrow {\simeq} X$ in which $F$ is a semi-flat $R$-complex.
\end{enumerate}
\end{definition}

Semi-projective, semi-injective, and semi-flat resolutions exist for any $R$-complex. Moreover, any right-bounded $R$-complex of projective (flat) modules is semi-projective (semi-flat), and any left-bounded $R$-complex of injective modules is semi-injective.

We now remind the total derived functors that we need.

\begin{remark} \label{3.3.2}
Let $\mathfrak{a}$ be an ideal of $R$, and $X$ and $Y$ two $R$-complexes. Then we have:
\begin{enumerate}
\item[(i)] Each of the functors $\Hom_{R}(X,-)$ and $\Hom_{R}(-,Y)$ on $\mathcal{C}(R)$ enjoys a right total derived functor on $\mathcal{D}(R)$, together with a balance property, in the sense that ${\bf R}\Hom_{R}(X,Y)$ can be computed by
$${\bf R}\Hom_{R}(X,Y)\simeq \Hom_{R}(P,Y) \simeq \Hom_{R}(X,I),$$
where $P\xrightarrow {\simeq} X$ is any semi-projective resolution of $X$, and $Y\xrightarrow {\simeq} I$ is any semi-injective resolution of $Y$. In addition, these functors turn out to be triangulated, in the sense that they preserve shifts and distinguished triangles. Moreover, we let $$\Ext^{i}_{R}(X,Y):=H_{-i}\left({\bf R}\Hom_{R}(X,Y)\right)$$
for every $i \in \mathbb{Z}$.
\item[(ii)] Each of the functors $X\otimes_{R}-$ and $-\otimes_{R}Y$ on $\mathcal{C}(R)$ enjoys a left total derived functor on $\mathcal{D}(R)$, together with a balance property, in the sense that $X\otimes_{R}^{\bf L}Y$ can be computed by
$$X\otimes_{R}^{\bf L}Y \simeq P\otimes_{R}Y \simeq X\otimes_{R}Q,$$
where $P\xrightarrow {\simeq} X$ is any semi-projective resolution of $X$, and $Q\xrightarrow {\simeq} Y$ is any semi-projective resolution of $Y$. Besides, these functors turn out to be triangulated. Moreover, we let $$\Tor^{R}_{i}(X,Y):=H_{i}\left(X\otimes_{R}^{\bf L}Y\right)$$
for every $i \in \mathbb{Z}$.
\item[(iii)] The functor $\Gamma_{\mathfrak{a}}(-)$ on $\mathcal{M}(R)$ extends naturally to a functor on $\mathcal{C}(R)$. The extended functor enjoys a right total derived functor ${\bf R}\Gamma_{\mathfrak{a}}(-):\mathcal{D}(R)\rightarrow \mathcal{D}(R)$, that can be computed by
${\bf R}\Gamma_{\mathfrak{a}}(X)\simeq \Gamma_{\mathfrak{a}}(I)$, where $X \xrightarrow {\simeq} I$ is any semi-injective resolution of $X$. Besides, we define the $i$th local cohomology module of $X$ to be
$$H^{i}_{\mathfrak{a}}(X):= H_{-i}\left({\bf R}\Gamma_{\mathfrak{a}}(X)\right)$$
for every $i\in\mathbb{Z}$. The functor ${\bf R}\Gamma_{\mathfrak{a}}(-)$ turns out to be triangulated.
\item[(iv)] The functor $\Lambda^{\mathfrak{a}}(-)$ on $\mathcal{M}(R)$ extends naturally to a functor on $\mathcal{C}(R)$. The extended functor enjoys a left total derived functor ${\bf L}\Lambda^{\mathfrak{a}}(-):\mathcal{D}(R)\rightarrow \mathcal{D}(R)$, that can be computed by ${\bf L}\Lambda^{\mathfrak{a}}(X)\simeq \Lambda^{\mathfrak{a}}(P)$, where $P \xrightarrow {\simeq} X$ is any semi-projective resolution of $X$. Moreover, we define the $i$th local homology module of $X$ to be
$$H^{\mathfrak{a}}_{i}(X):= H_{i}\left({\bf L}\Lambda^{\mathfrak{a}}(X)\right)$$
for every $i\in\mathbb{Z}$. The functor ${\bf L}\Lambda^{\mathfrak{a}}(-)$ turns out to be triangulated.
\end{enumerate}
\end{remark}

We further need the notion of way-out functors for functors between the category of complexes.

\begin{definition} \label{3.3.3}
Let $R$ and $S$ be two rings, and $\mathcal{F}: \mathcal{C}(R) \rightarrow \mathcal{C}(S)$ a covariant functor. Then:
\begin{enumerate}
\item[(i)] $\mathcal{F}$ is said to be \textit{way-out left} \index{way-out left functor on category of complexes} if for every $n \in \mathbb{Z}$, there is an $m \in \mathbb{Z}$, such that for any $R$-complex $X$ with $X_{i}=0$ for every $i>m$, we have $\mathcal{F}(X)_{i}=0$ for every $i>n$.
\item[(ii)] $\mathcal{F}$ is said to be \textit{way-out right} \index{way-out right functor on category of complexes} if for every $n \in \mathbb{Z}$, there is an $m \in \mathbb{Z}$, such that for any $R$-complex $X$ with $X_{i}=0$ for every $i<m$, we have $\mathcal{F}(X)_{i}=0$ for every $i<n$.
\item[(iii)] $\mathcal{F}$ is said to be \textit{way-out} \index{way-out functor on category of complexes} if it is both way-out left and way-out right.
\end{enumerate}
\end{definition}

The following lemma is the Way-out Lemma for functors between the category of complexes. We include a proof since there is no account of this version in the literature.

\begin{lemma} \label{3.3.4}
Let $R$ and $S$ be two rings, and $\mathcal{F},\mathcal{G}: \mathcal{C}(R) \rightarrow \mathcal{C}(S)$ two additive covariant functors that commute with shift and preserve the exactness of degreewise split short exact sequences of $R$-complexes. Let $\sigma : \mathcal{F} \rightarrow \mathcal{G}$ be a natural transformation of functors. Then the following assertions hold:
\begin{enumerate}
\item[(i)] If $X$ is a bounded $R$-complex such that $\sigma^{X_{i}}:\mathcal{F}(X_{i}) \rightarrow \mathcal{G}(X_{i})$ is a quasi-isomorphism for every $i \in \mathbb{Z}$, then $\sigma^{X}: \mathcal{F}(X) \rightarrow \mathcal{G}(X)$ is a quasi-isomorphism.
\item[(ii)] If $\mathcal{F}$ and $\mathcal{G}$ are way-out left, and $X$ is a left-bounded $R$-complex such that $\sigma^{X_{i}}:\mathcal{F}(X_{i}) \rightarrow \mathcal{G}(X_{i})$ is a quasi-isomorphism for every $i \in \mathbb{Z}$, then $\sigma^{X}: \mathcal{F}(X) \rightarrow \mathcal{G}(X)$ is a quasi-isomorphism.
\item[(iii)] If $\mathcal{F}$ and $\mathcal{G}$ are way-out right, and $X$ is a right-bounded $R$-complex such that $\sigma^{X_{i}}:\mathcal{F}(X_{i}) \rightarrow \mathcal{G}(X_{i})$ is a quasi-isomorphism for every $i \in \mathbb{Z}$, then $\sigma^{X}: \mathcal{F}(X) \rightarrow \mathcal{G}(X)$ is a quasi-isomorphism.
\item[(iv)] If $\mathcal{F}$ and $\mathcal{G}$ are way-out, and $X$ is an $R$-complex such that $\sigma^{X_{i}}:\mathcal{F}(X_{i}) \rightarrow \mathcal{G}(X_{i})$ is a quasi-isomorphism for every $i \in \mathbb{Z}$, then $\sigma^{X}: \mathcal{F}(X) \rightarrow \mathcal{G}(X)$ is a quasi-isomorphism.
\end{enumerate}
\end{lemma}

\begin{prf}
(i): Without loss of generality we may assume that
\[X: 0\rightarrow X_{n} \xrightarrow{\partial_{n}^{X}} X_{n-1} \rightarrow \cdots \rightarrow X_{1} \xrightarrow{\partial_{1}^{X}} X_{0} \rightarrow 0.\]
Let
\[Y: 0\rightarrow X_{n-1} \xrightarrow{\partial_{n-1}^{X}} X_{n-2} \rightarrow \cdots \rightarrow X_{1} \xrightarrow{\partial_{1}^{X}} X_{0} \rightarrow 0.\]
Consider the degreewise split short exact sequence
\[0\rightarrow Y \rightarrow X \rightarrow \Sigma^{n}X_{n} \rightarrow 0\]
of $R$-complexes, and apply $\mathcal{F}$ and $\mathcal{G}$ to get the following commutative diagram of $S$-complexes with exact rows:
\[\begin{tikzpicture}[every node/.style={midway},]
  \matrix[column sep={2.5em}, row sep={2.5em}]
  {\node(1) {$0$}; & \node(2) {$\mathcal{F}(Y)$}; & \node(3) {$\mathcal{F}(X)$}; & \node(4) {$\Sigma^{n}\mathcal{F}(X_{n})$}; & \node(5) {$0$};\\
  \node(6) {$0$}; & \node(7) {$\mathcal{G}(Y)$}; & \node(8) {$\mathcal{G}(X)$}; & \node(9) {$\Sigma^{n}\mathcal{G}(X_{n})$}; & \node(10) {$0$};\\};
  \draw[decoration={markings,mark=at position 1 with {\arrow[scale=1.5]{>}}},postaction={decorate},shorten >=0.5pt] (4) -- (9) node[anchor=west] {$\Sigma^{n}\sigma_{X_{n}}$};
  \draw[decoration={markings,mark=at position 1 with {\arrow[scale=1.5]{>}}},postaction={decorate},shorten >=0.5pt] (3) -- (8) node[anchor=west] {$\sigma_{X}$};
  \draw[decoration={markings,mark=at position 1 with {\arrow[scale=1.5]{>}}},postaction={decorate},shorten >=0.5pt] (2) -- (7) node[anchor=west] {$\sigma_{Y}$};
  \draw[decoration={markings,mark=at position 1 with {\arrow[scale=1.5]{>}}},postaction={decorate},shorten >=0.5pt] (1) -- (2) node[anchor=south] {};
  \draw[decoration={markings,mark=at position 1 with {\arrow[scale=1.5]{>}}},postaction={decorate},shorten >=0.5pt] (2) -- (3) node[anchor=south] {};
  \draw[decoration={markings,mark=at position 1 with {\arrow[scale=1.5]{>}}},postaction={decorate},shorten >=0.5pt] (3) -- (4) node[anchor=south] {};
  \draw[decoration={markings,mark=at position 1 with {\arrow[scale=1.5]{>}}},postaction={decorate},shorten >=0.5pt] (4) -- (5) node[anchor=south] {};
  \draw[decoration={markings,mark=at position 1 with {\arrow[scale=1.5]{>}}},postaction={decorate},shorten >=0.5pt] (6) -- (7) node[anchor=south] {};
  \draw[decoration={markings,mark=at position 1 with {\arrow[scale=1.5]{>}}},postaction={decorate},shorten >=0.5pt] (7) -- (8) node[anchor=south] {};
  \draw[decoration={markings,mark=at position 1 with {\arrow[scale=1.5]{>}}},postaction={decorate},shorten >=0.5pt] (8) -- (9) node[anchor=south] {};
  \draw[decoration={markings,mark=at position 1 with {\arrow[scale=1.5]{>}}},postaction={decorate},shorten >=0.5pt] (9) -- (10) node[anchor=south] {};
\end{tikzpicture}\]
Note that $\Sigma^{n}\sigma_{X_{n}}$ is a quasi-isomorphism by the assumption. Hence to prove that $\sigma_{X}$ is a quasi-isomorphism, it suffices to show that $\sigma_{Y}$ is a quasi-isomorphism. Since $Y$ is bounded, by continuing this process with $Y$, we reach at a level that we need $\sigma_{X_{0}}$ to be a quasi-isomorphism, which holds by the assumption. Therefore, we are done.

(ii): Without loss of generality we may assume that
\[X: 0\rightarrow X_{n} \xrightarrow{\partial_{n}^{X}} X_{n-1} \rightarrow \cdots.\]
Let $i \in \mathbb{Z}$. We show that $H_{i}(\sigma_{X}):H_{i}\left(\mathcal{F}(X)\right) \rightarrow H_{i}\left(\mathcal{G}(X)\right)$ is an isomorphism. Since $\mathcal{F}$ and $\mathcal{G}$ are way-out left, we can choose an integer $j \in \mathbb{Z}$ corresponding to $i-2$. Let
\[Z: 0\rightarrow X_{n} \xrightarrow{\partial^{X}_{n}} X_{n-1} \rightarrow \cdots \rightarrow X_{j+1} \xrightarrow{\partial^{X}_{j+1}} X_{j} \rightarrow 0\]
and
\[Y: 0 \rightarrow X_{j-1} \xrightarrow{\partial^{X}_{j-1}} X_{j-2} \rightarrow \cdots.\]
Then there is a degreewise split short exact sequence
\[0 \rightarrow Y \rightarrow X \rightarrow Z \rightarrow 0\]
of $R$-complexes. Apply $\mathcal{F}$ and $\mathcal{G}$ to get the following commutative diagram with exact rows:
\[\begin{tikzpicture}[every node/.style={midway},]
  \matrix[column sep={2.5em}, row sep={2.5em}]
  {\node(1) {$0$}; & \node(2) {$\mathcal{F}(Y)$}; & \node(3) {$\mathcal{F}(X)$}; & \node(4) {$\mathcal{F}(Z)$}; & \node(5) {$0$};\\
  \node(6) {$0$}; & \node(7) {$\mathcal{G}(Y)$}; & \node(8) {$\mathcal{G}(X)$}; & \node(9) {$\mathcal{G}(Z)$}; & \node(10) {$0$};\\};
  \draw[decoration={markings,mark=at position 1 with {\arrow[scale=1.5]{>}}},postaction={decorate},shorten >=0.5pt] (4) -- (9) node[anchor=west] {$\sigma_{Z}$};
  \draw[decoration={markings,mark=at position 1 with {\arrow[scale=1.5]{>}}},postaction={decorate},shorten >=0.5pt] (3) -- (8) node[anchor=west] {$\sigma_{X}$};
  \draw[decoration={markings,mark=at position 1 with {\arrow[scale=1.5]{>}}},postaction={decorate},shorten >=0.5pt] (2) -- (7) node[anchor=west] {$\sigma_{Y}$};
  \draw[decoration={markings,mark=at position 1 with {\arrow[scale=1.5]{>}}},postaction={decorate},shorten >=0.5pt] (1) -- (2) node[anchor=south] {};
  \draw[decoration={markings,mark=at position 1 with {\arrow[scale=1.5]{>}}},postaction={decorate},shorten >=0.5pt] (2) -- (3) node[anchor=south] {};
  \draw[decoration={markings,mark=at position 1 with {\arrow[scale=1.5]{>}}},postaction={decorate},shorten >=0.5pt] (3) -- (4) node[anchor=south] {};
  \draw[decoration={markings,mark=at position 1 with {\arrow[scale=1.5]{>}}},postaction={decorate},shorten >=0.5pt] (4) -- (5) node[anchor=south] {};
  \draw[decoration={markings,mark=at position 1 with {\arrow[scale=1.5]{>}}},postaction={decorate},shorten >=0.5pt] (6) -- (7) node[anchor=south] {};
  \draw[decoration={markings,mark=at position 1 with {\arrow[scale=1.5]{>}}},postaction={decorate},shorten >=0.5pt] (7) -- (8) node[anchor=south] {};
  \draw[decoration={markings,mark=at position 1 with {\arrow[scale=1.5]{>}}},postaction={decorate},shorten >=0.5pt] (8) -- (9) node[anchor=south] {};
  \draw[decoration={markings,mark=at position 1 with {\arrow[scale=1.5]{>}}},postaction={decorate},shorten >=0.5pt] (9) -- (10) node[anchor=south] {};
\end{tikzpicture}\]
From the above diagram, we get the following commutative diagram of $S$-modules with exact rows:
\[\begin{tikzpicture}[every node/.style={midway},]
  \matrix[column sep={2.5em}, row sep={2.5em}]
  {\node(1) {$0=H_{i}\left(\mathcal{F}(Y)\right)$}; & \node(2) {$H_{i}\left(\mathcal{F}(X)\right)$}; & \node(3) {$H_{i}\left(\mathcal{F}(Z)\right)$}; & \node(4) {$H_{i-1}\left(\mathcal{F}(Y)\right)=0$}; \\
  \node(5) {$0=H_{i}\left(\mathcal{G}(Y)\right)$}; & \node(6) {$H_{i}\left(\mathcal{G}(X)\right)$}; & \node(7) {$H_{i}\left(\mathcal{G}(Z)\right)$}; & \node(8) {$H_{i-1}\left(\mathcal{G}(Y)\right)=0$};\\};
  \draw[decoration={markings,mark=at position 1 with {\arrow[scale=1.5]{>}}},postaction={decorate},shorten >=0.5pt] (2) -- (6) node[anchor=west] {$H_{i}(\sigma_{X})$};
  \draw[decoration={markings,mark=at position 1 with {\arrow[scale=1.5]{>}}},postaction={decorate},shorten >=0.5pt] (3) -- (7) node[anchor=west] {$H_{i}(\sigma_{Z})$};
  \draw[decoration={markings,mark=at position 1 with {\arrow[scale=1.5]{>}}},postaction={decorate},shorten >=0.5pt] (1) -- (2) node[anchor=south] {};
  \draw[decoration={markings,mark=at position 1 with {\arrow[scale=1.5]{>}}},postaction={decorate},shorten >=0.5pt] (2) -- (3) node[anchor=south] {};
  \draw[decoration={markings,mark=at position 1 with {\arrow[scale=1.5]{>}}},postaction={decorate},shorten >=0.5pt] (3) -- (4) node[anchor=south] {};
  \draw[decoration={markings,mark=at position 1 with {\arrow[scale=1.5]{>}}},postaction={decorate},shorten >=0.5pt] (5) -- (6) node[anchor=south] {};
  \draw[decoration={markings,mark=at position 1 with {\arrow[scale=1.5]{>}}},postaction={decorate},shorten >=0.5pt] (6) -- (7) node[anchor=south] {};
  \draw[decoration={markings,mark=at position 1 with {\arrow[scale=1.5]{>}}},postaction={decorate},shorten >=0.5pt] (7) -- (8) node[anchor=south] {};
\end{tikzpicture}\]
where the vanishing is due to the choice of $j$. Since $Z$ is bounded, it follows from (i) that $H_{i}(\sigma_{Z})$ is an isomorphism, and as a consequence, $H_{i}(\sigma_{X})$ is an isomorphism.

(iii): Without loss of generality we may assume that
\[X: \cdots \rightarrow X_{n+1} \xrightarrow{\partial_{n+1}^{X}} X_{n} \rightarrow 0.\]
Let $i \in \mathbb{Z}$. We show that $H_{i}(\sigma_{X}):H_{i}\left(\mathcal{F}(X)\right) \rightarrow H_{i}\left(\mathcal{G}(X)\right)$ is an isomorphism. Since $\mathcal{F}$ and $\mathcal{G}$ are way-out right, we can choose an integer $j \in \mathbb{Z}$ corresponding to $i+2$. Let
\[Y: 0\rightarrow X_{j-1} \xrightarrow{\partial^{X}_{j-1}} X_{j-2} \rightarrow \cdots \rightarrow X_{n+1} \xrightarrow{\partial^{X}_{n+1}} X_{n} \rightarrow 0\]
and
\[Z: \cdots \rightarrow X_{j+1} \xrightarrow{\partial^{X}_{j+1}} X_{j} \rightarrow 0.\]
Then there is a degreewise split short exact sequence
\[0 \rightarrow Y \rightarrow X \rightarrow Z \rightarrow 0\]
of $R$-complexes. Apply $\mathcal{F}$ and $\mathcal{G}$ to get the following commutative diagram of $S$-complexes with exact rows:
\[\begin{tikzpicture}[every node/.style={midway},]
  \matrix[column sep={2.5em}, row sep={2.5em}]
  {\node(1) {$0$}; & \node(2) {$\mathcal{F}(Y)$}; & \node(3) {$\mathcal{F}(X)$}; & \node(4) {$\mathcal{F}(Z)$}; & \node(5) {$0$};\\
  \node(6) {$0$}; & \node(7) {$\mathcal{G}(Y)$}; & \node(8) {$\mathcal{G}(X)$}; & \node(9) {$\mathcal{G}(Z)$}; & \node(10) {$0$};\\};
  \draw[decoration={markings,mark=at position 1 with {\arrow[scale=1.5]{>}}},postaction={decorate},shorten >=0.5pt] (4) -- (9) node[anchor=west] {$\sigma_{Z}$};
  \draw[decoration={markings,mark=at position 1 with {\arrow[scale=1.5]{>}}},postaction={decorate},shorten >=0.5pt] (3) -- (8) node[anchor=west] {$\sigma_{X}$};
  \draw[decoration={markings,mark=at position 1 with {\arrow[scale=1.5]{>}}},postaction={decorate},shorten >=0.5pt] (2) -- (7) node[anchor=west] {$\sigma_{Y}$};
  \draw[decoration={markings,mark=at position 1 with {\arrow[scale=1.5]{>}}},postaction={decorate},shorten >=0.5pt] (1) -- (2) node[anchor=south] {};
  \draw[decoration={markings,mark=at position 1 with {\arrow[scale=1.5]{>}}},postaction={decorate},shorten >=0.5pt] (2) -- (3) node[anchor=south] {};
  \draw[decoration={markings,mark=at position 1 with {\arrow[scale=1.5]{>}}},postaction={decorate},shorten >=0.5pt] (3) -- (4) node[anchor=south] {};
  \draw[decoration={markings,mark=at position 1 with {\arrow[scale=1.5]{>}}},postaction={decorate},shorten >=0.5pt] (4) -- (5) node[anchor=south] {};
  \draw[decoration={markings,mark=at position 1 with {\arrow[scale=1.5]{>}}},postaction={decorate},shorten >=0.5pt] (6) -- (7) node[anchor=south] {};
  \draw[decoration={markings,mark=at position 1 with {\arrow[scale=1.5]{>}}},postaction={decorate},shorten >=0.5pt] (7) -- (8) node[anchor=south] {};
  \draw[decoration={markings,mark=at position 1 with {\arrow[scale=1.5]{>}}},postaction={decorate},shorten >=0.5pt] (8) -- (9) node[anchor=south] {};
  \draw[decoration={markings,mark=at position 1 with {\arrow[scale=1.5]{>}}},postaction={decorate},shorten >=0.5pt] (9) -- (10) node[anchor=south] {};
\end{tikzpicture}\]
From the above diagram, we get the following commutative diagram of $S$-modules with exact rows:
\[\begin{tikzpicture}[every node/.style={midway},]
  \matrix[column sep={2.5em}, row sep={2.5em}]
  {\node(1) {$0=H_{i+1}\left(\mathcal{F}(Z)\right)$}; & \node(2) {$H_{i}\left(\mathcal{F}(Y)\right)$}; & \node(3) {$H_{i}\left(\mathcal{F}(X)\right)$}; & \node(4) {$H_{i}\left(\mathcal{F}(Z)\right)=0$}; \\
  \node(5) {$0=H_{i+1}\left(\mathcal{G}(Z)\right)$}; & \node(6) {$H_{i}\left(\mathcal{G}(Y)\right)$}; & \node(7) {$H_{i}\left(\mathcal{G}(X)\right)$}; & \node(8) {$H_{i}\left(\mathcal{G}(Z)\right)=0$};\\};
  \draw[decoration={markings,mark=at position 1 with {\arrow[scale=1.5]{>}}},postaction={decorate},shorten >=0.5pt] (2) -- (6) node[anchor=west] {$H_{i}(\sigma_{Y})$};
  \draw[decoration={markings,mark=at position 1 with {\arrow[scale=1.5]{>}}},postaction={decorate},shorten >=0.5pt] (3) -- (7) node[anchor=west] {$H_{i}(\sigma_{X})$};
  \draw[decoration={markings,mark=at position 1 with {\arrow[scale=1.5]{>}}},postaction={decorate},shorten >=0.5pt] (1) -- (2) node[anchor=south] {};
  \draw[decoration={markings,mark=at position 1 with {\arrow[scale=1.5]{>}}},postaction={decorate},shorten >=0.5pt] (2) -- (3) node[anchor=south] {};
  \draw[decoration={markings,mark=at position 1 with {\arrow[scale=1.5]{>}}},postaction={decorate},shorten >=0.5pt] (3) -- (4) node[anchor=south] {};
  \draw[decoration={markings,mark=at position 1 with {\arrow[scale=1.5]{>}}},postaction={decorate},shorten >=0.5pt] (5) -- (6) node[anchor=south] {};
  \draw[decoration={markings,mark=at position 1 with {\arrow[scale=1.5]{>}}},postaction={decorate},shorten >=0.5pt] (6) -- (7) node[anchor=south] {};
  \draw[decoration={markings,mark=at position 1 with {\arrow[scale=1.5]{>}}},postaction={decorate},shorten >=0.5pt] (7) -- (8) node[anchor=south] {};
\end{tikzpicture}\]
where the vanishing is due to the choice of $j$. Since $Y$ is bounded, it follows from (i) that $H_{i}(\sigma_{Y})$ is an isomorphism, and as a consequence, $H_{i}(\sigma_{X})$ is an isomorphism.

(iv): Let
\[Y: 0 \rightarrow X_{0} \xrightarrow{\partial^{X}_{0}} X_{-1} \rightarrow \cdots\]
and
\[Z: \cdots \rightarrow X_{2} \xrightarrow{\partial^{X}_{2}} X_{1} \rightarrow 0.\]
Then there is a degreewise split short exact sequence
\[0 \rightarrow Y \rightarrow X \rightarrow Z \rightarrow 0\]
of $R$-complexes.
Applying $\mathcal{F}$ and $\mathcal{G}$, we get the following commutative diagram of $S$-complexes with exact rows:
\[\begin{tikzpicture}[every node/.style={midway},]
  \matrix[column sep={2.5em}, row sep={2.5em}]
  {\node(1) {$0$}; & \node(2) {$\mathcal{F}(Y)$}; & \node(3) {$\mathcal{F}(X)$}; & \node(4) {$\mathcal{F}(Z)$}; & \node(5) {$0$};\\
  \node(6) {$0$}; & \node(7) {$\mathcal{G}(Y)$}; & \node(8) {$\mathcal{G}(X)$}; & \node(9) {$\mathcal{G}(Z)$}; & \node(10) {$0$};\\};
  \draw[decoration={markings,mark=at position 1 with {\arrow[scale=1.5]{>}}},postaction={decorate},shorten >=0.5pt] (4) -- (9) node[anchor=west] {$\sigma_{Z}$};
  \draw[decoration={markings,mark=at position 1 with {\arrow[scale=1.5]{>}}},postaction={decorate},shorten >=0.5pt] (3) -- (8) node[anchor=west] {$\sigma_{X}$};
  \draw[decoration={markings,mark=at position 1 with {\arrow[scale=1.5]{>}}},postaction={decorate},shorten >=0.5pt] (2) -- (7) node[anchor=west] {$\sigma_{Y}$};
  \draw[decoration={markings,mark=at position 1 with {\arrow[scale=1.5]{>}}},postaction={decorate},shorten >=0.5pt] (1) -- (2) node[anchor=south] {};
  \draw[decoration={markings,mark=at position 1 with {\arrow[scale=1.5]{>}}},postaction={decorate},shorten >=0.5pt] (2) -- (3) node[anchor=south] {};
  \draw[decoration={markings,mark=at position 1 with {\arrow[scale=1.5]{>}}},postaction={decorate},shorten >=0.5pt] (3) -- (4) node[anchor=south] {};
  \draw[decoration={markings,mark=at position 1 with {\arrow[scale=1.5]{>}}},postaction={decorate},shorten >=0.5pt] (4) -- (5) node[anchor=south] {};
  \draw[decoration={markings,mark=at position 1 with {\arrow[scale=1.5]{>}}},postaction={decorate},shorten >=0.5pt] (6) -- (7) node[anchor=south] {};
  \draw[decoration={markings,mark=at position 1 with {\arrow[scale=1.5]{>}}},postaction={decorate},shorten >=0.5pt] (7) -- (8) node[anchor=south] {};
  \draw[decoration={markings,mark=at position 1 with {\arrow[scale=1.5]{>}}},postaction={decorate},shorten >=0.5pt] (8) -- (9) node[anchor=south] {};
  \draw[decoration={markings,mark=at position 1 with {\arrow[scale=1.5]{>}}},postaction={decorate},shorten >=0.5pt] (9) -- (10) node[anchor=south] {};
\end{tikzpicture}\]
Since $Y$ is left-bounded, $\sigma_{Y}$ is a quasi-isomorphism by (ii), and since $Z$ is right-bounded, $\sigma_{Z}$ is a quasi-isomorphism by (iii). Therefore, $\sigma_{X}$ is a quasi-isomorphism.
\end{prf}

Although $\check{C}_{\infty}(\underline{a})$ is suitable in Proposition \ref{3.2.12}, it is not applicable in the next proposition due to the fact that it is concentrated in degrees $1,0,...,-n$. What we really need here is a semi-projective approximation of $\check{C}(\underline{a})$ of the same length, i.e. concentrated in degrees $0,-1,...,-n$. We proceed as follows.

Given an element $a \in R$, consider the following commutative diagram:
\[\begin{tikzpicture}[every node/.style={midway},]
  \matrix[column sep={3em}, row sep={3em}]
  {\node(1) {$0$}; & \node(2) {$R[X]\oplus R$}; & \node(3) {$R[X]$}; & \node(4) {$0$};\\
  \node(5) {$0$}; & \node(6) {$R$}; & \node(7) {$R_{a}$}; & \node(8) {$0$};\\};
  \draw[decoration={markings,mark=at position 1 with {\arrow[scale=1.5]{>}}},postaction={decorate},shorten >=0.5pt] (1) -- (2) node[anchor=east] {};
  \draw[decoration={markings,mark=at position 1 with {\arrow[scale=1.5]{>}}},postaction={decorate},shorten >=0.5pt] (2) -- (3) node[anchor=south] {$f_{a}$};
  \draw[decoration={markings,mark=at position 1 with {\arrow[scale=1.5]{>}}},postaction={decorate},shorten >=0.5pt] (3) -- (4) node[anchor=east] {};
  \draw[decoration={markings,mark=at position 1 with {\arrow[scale=1.5]{>}}},postaction={decorate},shorten >=0.5pt] (5) -- (6) node[anchor=south] {};
  \draw[decoration={markings,mark=at position 1 with {\arrow[scale=1.5]{>}}},postaction={decorate},shorten >=0.5pt] (6) -- (7) node[anchor=south] {$\lambda_{R}^{a}$};
  \draw[decoration={markings,mark=at position 1 with {\arrow[scale=1.5]{>}}},postaction={decorate},shorten >=0.5pt] (7) -- (8) node[anchor=south] {};
  \draw[decoration={markings,mark=at position 1 with {\arrow[scale=1.5]{>}}},postaction={decorate},shorten >=0.5pt] (2) -- (6) node[anchor=west] {$\pi$};
  \draw[decoration={markings,mark=at position 1 with {\arrow[scale=1.5]{>}}},postaction={decorate},shorten >=0.5pt] (3) -- (7) node[anchor=west] {$g_{a}$};
\end{tikzpicture}\]
in which, $f_{a}\left(p(X),b\right)=(aX-1)p(X)+b$, $\pi\left(p(X),b\right)=b$, $\lambda^{a}_{R}$ is the localization map, and $g_{a}\left(p(X)\right)=\frac{b_{k}}{a^{k}}+\cdots+\frac{b_{1}}{a}+\frac{b_{0}}{1}$ where $p(X)=b_{k}X^{k}+\cdots+b_{1}X+b_{0} \in R[X]$. Let $L^{R}(a)$ denote the $R$-complex in the first row of the diagram above concentrated in degrees $0,-1$. Since the second row is isomorphic to $\check{C}(a)$, it can be seen that the diagram above provides a quasi-isomorphism $L^{R}(a) \xrightarrow{\simeq} \check{C}(a)$. Hence $L^{R}(a) \xrightarrow{\simeq} \check{C}(a)$ is a semi-projective resolution of $\check{C}(a)$. Now for the elements $\underline{�}=a_{1},...,a_{n}\in R$, let
$$L^{R}(\underline{a})=L^{R}(a_{1})\otimes_{R}\cdots\otimes_{R}L^{R}(a_{n}).$$
Then $L^{R}(\underline{a})$ is an $R$-complex of free modules concentrated in degrees $0,-1,...,-n$, and $L^{R}(\underline{a}) \xrightarrow{\simeq} \check{C}(\underline{a})$ is a semi-projective resolution of $\check{C}(\underline{a})$.

The next proposition inspects the relation between derived torsion functor and derived completion functor with \v{C}ech complex, and provides the second crucial step towards the Greenlees-May Duality.

\begin{proposition} \label{3.3.5}
Let $\mathfrak{a}= (a_{1},...,a_{n})$ be an ideal of $R$, $\underline{a}=a_{1},...,a_{n}$, and $X \in \mathcal{D}(R)$. Then there are natural isomorphisms in $\mathcal{D}(R)$:
\begin{enumerate}
\item[(i)] ${\bf R}\Gamma_{\mathfrak{a}}(X) \simeq \check{C}(\underline{a})\otimes_{R}^{\bf L}X \simeq \check{C}_{\infty}(\underline{a})\otimes_{R}^{\bf L}X$.
\item[(ii)] ${\bf L}\Lambda^{\mathfrak{a}}(X) \simeq {\bf R}\Hom_{R}\left(\check{C}(\underline{a}),X\right) \simeq {\bf R}\Hom_{R}\left(\check{C}_{\infty}(\underline{a}),X\right)$.
\end{enumerate}
\end{proposition}

\begin{prf}
(i): Let $X \xrightarrow{\simeq} I$ be a semi-injective resolution of $X$. Then ${\bf R}\Gamma_{\mathfrak{a}}(X) \simeq \Gamma_{\mathfrak{a}}(I)$, and
$$\check{C}(\underline{a})\otimes_{R}^{\bf L}X \simeq \check{C}(\underline{a})\otimes_{R}^{\bf L}I \simeq \check{C}(\underline{a})\otimes_{R}I,$$
since $\check{C}(\underline{a})$ is a semi-flat $R$-complex. Hence it suffices to establish a quasi-isomorphism $\Gamma_{\mathfrak{a}}(I) \rightarrow \check{C}(\underline{a})\otimes_{R}I$.

Let $Y$ be an $R$-complex and $i \in \mathbb{Z}$. Let $\sigma_{i}^{Y}: \Gamma_{\mathfrak{a}}(Y)_{i} \rightarrow \left(\check{C}(\underline{a})\otimes_{R}Y\right)_{i}$ be the composition of the following natural $R$-homomorphisms:
$$\Gamma_{\mathfrak{a}}(Y)_{i} = \Gamma_{\mathfrak{a}}(Y_{i}) \xrightarrow{\cong} H^{0}_{\mathfrak{a}}(Y_{i}) \xrightarrow{\cong} H_{0}\left(\check{C}(\underline{a})\otimes_{R}Y_{i}\right)= \ker \left(\partial_{0}^{\check{C}(\underline{a})}\otimes_{R}Y_{i}\right) \rightarrow \check{C}(\underline{a})_{0}\otimes_{R}Y_{i} $$$$ \rightarrow \bigoplus_{s+t=i}\left(\check{C}(\underline{a})_{s}\otimes_{R}Y_{t}\right) = \left(\check{C}(\underline{a})\otimes_{R}Y\right)_{i}$$
We note that the second isomorphism above comes from Proposition \ref{3.2.12} (i). One can easily see that $\sigma^{Y}=(\sigma_{i}^{Y})_{i \in \mathbb{Z}}:\Gamma_{\mathfrak{a}}(Y) \rightarrow \check{C}(\underline{a})\otimes_{R}Y$ is a natural morphism of $R$-complexes.

Since $I_{i}$ is an injective $R$-module for any $i \in \mathbb{Z}$, using Proposition \ref{3.2.12} (i), we get
$$H_{-j}\left(\check{C}(\underline{a})\otimes_{R}I_{i}\right) \cong H^{j}_{\mathfrak{a}}(I_{i})=0$$
for every $j \geq 1$. It follows that $\sigma^{I_{i}}:\Gamma_{\mathfrak{a}}(I_{i}) \rightarrow \check{C}(\underline{a})\otimes_{R}I_{i}$ is a quasi-isomorphism:
\[\begin{tikzpicture}[every node/.style={midway}]
  \matrix[column sep={3em}, row sep={3em}]
  {\node(1) {$0$}; & \node(2) {$\Gamma_{\mathfrak{a}}(I_{i})$}; & \node(3) {$0$}; & \node(4) {$\cdots$}; & \node(6) {$0$}; & \node(7) {$0$};\\
  \node(8) {$0$}; & \node(9) {$\check{C}(\underline{a})_{0}\otimes_{R}I_{i}$}; & \node(10) {$\check{C}(\underline{a})_{-1}\otimes_{R}I_{i}$}; & \node(11) {$\cdots$}; & \node(13) {$\check{C}(\underline{a})_{-n}\otimes_{R}I_{i}$}; & \node(14) {$0$};\\};
  \draw[decoration={markings,mark=at position 1 with {\arrow[scale=1.5]{>}}},postaction={decorate},shorten >=0.5pt] (2) -- (9) node[anchor=west] {$\sigma^{I_{i}}_{0}$};
  \draw[decoration={markings,mark=at position 1 with {\arrow[scale=1.5]{>}}},postaction={decorate},shorten >=0.5pt] (3) -- (10) node[anchor=west] {};
  \draw[decoration={markings,mark=at position 1 with {\arrow[scale=1.5]{>}}},postaction={decorate},shorten >=0.5pt] (6) -- (13) node[anchor=west] {};
  \draw[decoration={markings,mark=at position 1 with {\arrow[scale=1.5]{>}}},postaction={decorate},shorten >=0.5pt] (1) -- (2) node[anchor=south] {};
  \draw[decoration={markings,mark=at position 1 with {\arrow[scale=1.5]{>}}},postaction={decorate},shorten >=0.5pt] (8) -- (9) node[anchor=south] {};
  \draw[decoration={markings,mark=at position 1 with {\arrow[scale=1.5]{>}}},postaction={decorate},shorten >=0.5pt] (2) -- (3) node[anchor=south] {};
  \draw[decoration={markings,mark=at position 1 with {\arrow[scale=1.5]{>}}},postaction={decorate},shorten >=0.5pt] (9) -- (10) node[anchor=south] {};
  \draw[decoration={markings,mark=at position 1 with {\arrow[scale=1.5]{>}}},postaction={decorate},shorten >=0.5pt] (3) -- (4) node[anchor=south] {};
  \draw[decoration={markings,mark=at position 1 with {\arrow[scale=1.5]{>}}},postaction={decorate},shorten >=0.5pt] (10) -- (11) node[anchor=south] {};
  \draw[decoration={markings,mark=at position 1 with {\arrow[scale=1.5]{>}}},postaction={decorate},shorten >=0.5pt] (4) -- (6) node[anchor=south] {};
  \draw[decoration={markings,mark=at position 1 with {\arrow[scale=1.5]{>}}},postaction={decorate},shorten >=0.5pt] (11) -- (13) node[anchor=south] {};
  \draw[decoration={markings,mark=at position 1 with {\arrow[scale=1.5]{>}}},postaction={decorate},shorten >=0.5pt] (6) -- (7) node[anchor=south] {};
  \draw[decoration={markings,mark=at position 1 with {\arrow[scale=1.5]{>}}},postaction={decorate},shorten >=0.5pt] (13) -- (14) node[anchor=south] {};
\end{tikzpicture}\]
In addition, it is easily seen that the functors $\Gamma_{\mathfrak{a}}(-):\mathcal{C}(R) \rightarrow \mathcal{C}(R)$ and $\check{C}(\underline{a})\otimes_{R}-:\mathcal{C}(R) \rightarrow \mathcal{C}(R)$ are additive way-out functors that commute with shift and preserve the exactness of degreewise split short exact sequences of $R$-complexes. Hence by Lemma \ref{3.3.4} (iv), we conclude that $\sigma^{I}:\Gamma_{\mathfrak{a}}(I) \rightarrow \check{C}(\underline{a})\otimes_{R}I$ is a quasi-isomorphism.

The second isomorphism is immediate since $\check{C}(\underline{a}) \simeq \check{C}_{\infty}(\underline{a})$ and $-\otimes_{R}^{\bf L}X$ is a functor on $\mathcal{D}(R)$.

(ii): We know that $L^{R}(\underline{a}) \simeq \check{C}(\underline{a}) \simeq \check{C}_{\infty}(\underline{a})$. Let $P \xrightarrow{\simeq} X$ be a semi-projective resolution of $X$. Then ${\bf L}\Lambda^{\mathfrak{a}}(X) \simeq \Lambda^{\mathfrak{a}}(P)$, and
$${\bf R}\Hom_{R}\left(\check{C}(\underline{a}),X\right) \simeq {\bf R}\Hom_{R}\left(L^{R}(\underline{a}),P\right) \simeq \Hom_{R}\left(L^{R}(\underline{a}),P\right),$$
since $L^{R}(\underline{a})$ is a semi-projective $R$-complex. Moreover, we have
$$\Hom_{R}\left(L^{R}(\underline{a}),P\right) \simeq {\bf R}\Hom_{R}\left(L^{R}(\underline{a}),P\right) \simeq {\bf R}\Hom_{R}\left(\check{C}_{\infty}(\underline{a}),P\right)\simeq \Hom_{R}\left(\check{C}_{\infty}(\underline{a}),P\right),$$
since $\check{C}_{\infty}(\underline{a})$ is a semi-projective $R$-complex. In particular, we get
\begin{equation} \label{eq:3.3.5.1}
H_{i}\left(\Hom_{R}\left(L^{R}(\underline{a}),P\right)\right) \cong H_{i}\left(\Hom_{R}\left(\check{C}_{\infty}(\underline{a}),P\right)\right)
\end{equation}
for every $i \in \mathbb{Z}$. Now it suffices to establish a natural quasi-isomorphism $\Hom_{R}\left(L^{R}(\underline{a}),P\right) \rightarrow \Lambda^{\mathfrak{a}}(P)$.

Let $Y$ be an $R$-complex and $i \in \mathbb{Z}$. Let $\varsigma^{Y}_{i}: \Hom_{R}\left(L^{R}(\underline{a}),Y\right)_{i} \rightarrow \Lambda^{\mathfrak{a}}(Y)_{i}$ be the composition of the following natural $R$-homomorphisms:
$$\Hom_{R}\left(L^{R}(\underline{a}),Y\right)_{i} = \prod_{s \in \mathbb{Z}}\Hom_{R}\left(L^{R}(\underline{a})_{s},Y_{s+i}\right) \rightarrow \Hom_{R}\left(L^{R}(\underline{a})_{0},Y_{i}\right) $$$$ \rightarrow \frac{\Hom_{R}\left(L^{R}(\underline{a})_{0},Y_{i}\right)}{\im \left(\Hom_{R}\left(\partial_{0}^{L^{R}(\underline{a})},Y_{i}\right)\right)} = H_{0}\left(\Hom_{R}\left(L^{R}(\underline{a}),Y_{i}\right)\right) \xrightarrow{\cong} H_{0}\left(\Hom_{R}\left(\check{C}_{\infty}(\underline{a}),Y_{i}\right)\right) $$$$ \xrightarrow{\cong} H^{\mathfrak{a}}_{0}(Y_{i}) \rightarrow \Lambda^{\mathfrak{a}}(Y_{i}) = \Lambda^{\mathfrak{a}}(Y)_{i}$$
We note that the first isomorphism above comes from the isomorphism \eqref{eq:3.3.5.1} and the second comes from Proposition \ref{3.2.12} (ii). One can easily see that
$\varsigma^{Y}=(\varsigma^{Y}_{i})_{i \in \mathbb{Z}} : \Hom_{R}\left(L^{R}(\underline{a}),Y\right)\rightarrow \Lambda^{\mathfrak{a}}(Y)$ is a natural morphism of $R$-complexes.

Since $P_{i}$ is a projective $R$-module for any $i \in \mathbb{Z}$, using the isomorphism \eqref{eq:3.3.5.1} and Proposition \ref{3.2.12} (ii), we get
$$H_{j}\left(\Hom_{R}\left(L^{R}(\underline{a}),P_{i}\right)\right) \cong H_{j}\left(\Hom_{R}\left(\check{C}_{\infty}(\underline{a}),P_{i}\right)\right) \cong H^{\mathfrak{a}}_{j}(P_{i})=0$$
for every $j \geq 1$. It follows that $\varsigma^{P_{i}}: \Hom_{R}\left(L^{R}(\underline{a}),P_{i}\right) \rightarrow \Lambda^{\mathfrak{a}}(P_{i})$ is a quasi-isomorphism:
\[\begin{tikzpicture}[every node/.style={midway}]
  \matrix[column sep={1em}, row sep={3em}]
  {\node(1) {$0$}; & \node(2) {$\Hom_{R}\left(L^{R}(\underline{a})_{-n},P_{i}\right)$}; & \node(3) {$\cdots$}; & \node(4) {$\Hom_{R}\left(L^{R}(\underline{a})_{-1},P_{i}\right)$}; & \node(6) {$\Hom_{R}\left(L^{R}(\underline{a})_{0},P_{i}\right)$}; & \node(7) {$0$};\\
  \node(8) {$0$}; & \node(9) {$0$}; & \node(10) {$\cdots$}; & \node(11) {$0$}; & \node(13) {$\Lambda^{\mathfrak{a}}(P_{i})$}; & \node(14) {$0$};\\};
  \draw[decoration={markings,mark=at position 1 with {\arrow[scale=1.5]{>}}},postaction={decorate},shorten >=0.5pt] (2) -- (9) node[anchor=west] {};
  \draw[decoration={markings,mark=at position 1 with {\arrow[scale=1.5]{>}}},postaction={decorate},shorten >=0.5pt] (4) -- (11) node[anchor=west] {};
  \draw[decoration={markings,mark=at position 1 with {\arrow[scale=1.5]{>}}},postaction={decorate},shorten >=0.5pt] (6) -- (13) node[anchor=west] {$\varsigma^{P_{i}}_{0}$};
  \draw[decoration={markings,mark=at position 1 with {\arrow[scale=1.5]{>}}},postaction={decorate},shorten >=0.5pt] (1) -- (2) node[anchor=south] {};
  \draw[decoration={markings,mark=at position 1 with {\arrow[scale=1.5]{>}}},postaction={decorate},shorten >=0.5pt] (8) -- (9) node[anchor=south] {};
  \draw[decoration={markings,mark=at position 1 with {\arrow[scale=1.5]{>}}},postaction={decorate},shorten >=0.5pt] (2) -- (3) node[anchor=south] {};
  \draw[decoration={markings,mark=at position 1 with {\arrow[scale=1.5]{>}}},postaction={decorate},shorten >=0.5pt] (9) -- (10) node[anchor=south] {};
  \draw[decoration={markings,mark=at position 1 with {\arrow[scale=1.5]{>}}},postaction={decorate},shorten >=0.5pt] (3) -- (4) node[anchor=south] {};
  \draw[decoration={markings,mark=at position 1 with {\arrow[scale=1.5]{>}}},postaction={decorate},shorten >=0.5pt] (10) -- (11) node[anchor=south] {};
  \draw[decoration={markings,mark=at position 1 with {\arrow[scale=1.5]{>}}},postaction={decorate},shorten >=0.5pt] (4) -- (6) node[anchor=south] {};
  \draw[decoration={markings,mark=at position 1 with {\arrow[scale=1.5]{>}}},postaction={decorate},shorten >=0.5pt] (11) -- (13) node[anchor=south] {};
  \draw[decoration={markings,mark=at position 1 with {\arrow[scale=1.5]{>}}},postaction={decorate},shorten >=0.5pt] (6) -- (7) node[anchor=south] {};
  \draw[decoration={markings,mark=at position 1 with {\arrow[scale=1.5]{>}}},postaction={decorate},shorten >=0.5pt] (13) -- (14) node[anchor=south] {};
\end{tikzpicture}\]
In addition, it is easily seen that the functors $\Hom_{R}\left(L^{R}(\underline{a}),-\right):\mathcal{C}(R) \rightarrow \mathcal{C}(R)$ and $\Lambda^{\mathfrak{a}}(-):\mathcal{C}(R) \rightarrow \mathcal{C}(R)$ are additive way-out functors that commute with shift and preserve the exactness of degreewise split short exact sequences of $R$-complexes. Hence by Lemma \ref{3.3.4} (iv), we conclude that $\varsigma^{P}:\Hom_{R}\left(L^{R}(\underline{a}),P\right) \rightarrow \Lambda^{\mathfrak{a}}(P)$ is a quasi-isomorphism.

The second isomorphism is immediate since $\check{C}(\underline{a}) \simeq \check{C}_{\infty}(\underline{a})$ and ${\bf R}\Hom_{R}(-,X)$ is a functor on $\mathcal{D}(R)$.
\end{prf}

We note that if $\mathfrak{a}= (a_{1},...,a_{n})$ is an ideal of $R$ and $\underline{a}=a_{1},...,a_{n}$, then $\check{C}(\underline{a})$ as an element of $\mathcal{C}(R)$ depends on the generators $\underline{a}$. However, the proof of the next corollary shows that $\check{C}(\underline{a})$ as an element of $\mathcal{D}(R)$ is independent of the generators $\underline{a}$.

\begin{corollary} \label{3.3.6}
Let $\mathfrak{a}$ be an ideal of $R$. Then there are natural isomorphisms in $\mathcal{D}(R)$:
\begin{enumerate}
\item[(i)] ${\bf R}\Gamma_{\mathfrak{a}}(X) \simeq {\bf R}\Gamma_{\mathfrak{a}}(R)\otimes_{R}^{\bf L}X$.
\item[(ii)] ${\bf L}\Lambda^{\mathfrak{a}}(X) \simeq {\bf R}\Hom_{R}\left({\bf R}\Gamma_{\mathfrak{a}}(R),X\right)$.
\end{enumerate}
\end{corollary}

\begin{prf}
Suppose that $\mathfrak{a}= (a_{1},...,a_{n})$, and $\underline{a}=a_{1},...,a_{n}$. By Proposition \ref{3.3.5} (i), we have
$${\bf R}\Gamma_{\mathfrak{a}}(R) \simeq \check{C}(\underline{a})\otimes_{R}^{\bf L}R \simeq \check{C}(\underline{a}).$$
Now (i) and (ii) follow from Proposition \ref{3.3.5}.
\end{prf}

\section{Greenlees-May Duality}

Having the material developed in Sections 2 and 3 at our disposal, we are fully prepared to prove the celebrated Greenlees-May Duality Theorem. \index{Greenlees-May duality theorem}

\begin{theorem} \label{3.4.1}
Let $\mathfrak{a}$ be an ideal of $R$, and $X,Y \in \mathcal{D}(R)$. Then there is a natural isomorphism
$${\bf R}\Hom_{R}\left({\bf R}\Gamma_{\mathfrak{a}}(X),Y\right) \simeq {\bf R}\Hom_{R}\left(X,{\bf L}\Lambda^{\mathfrak{a}}(Y)\right)$$
in $\mathcal{D}(R)$.
\end{theorem}

\begin{prf}
Using Corollary \ref{3.3.6} and the Adjointness Isomorphism, we have
\begin{equation*}
\begin{split}
{\bf R}\Hom_{R}\left({\bf R}\Gamma_{\mathfrak{a}}(X),Y\right) & \simeq {\bf R}\Hom_{R}\left({\bf R}\Gamma_{\mathfrak{a}}(R)\otimes_{R}^{\bf L}X,Y\right) \\
 & \simeq {\bf R}\Hom_{R}\left(X,{\bf R}\Hom_{R}\left({\bf R}\Gamma_{\mathfrak{a}}(R),X\right)\right) \\
 & \simeq {\bf R}\Hom_{R}\left(X,{\bf L}\Lambda^{\mathfrak{a}}(Y)\right).
\end{split}
\end{equation*}
\end{prf}

\begin{corollary}
Let $\mathfrak{a}$ be an ideal of $R$, and $X,Y \in \mathcal{D}(R)$. Then there are natural isomorphisms:
\begin{equation*}
\begin{split}
{\bf L}\Lambda^{\mathfrak{a}}\left({\bf R}\Hom_{R}(X,Y)\right) & \simeq {\bf R}\Hom_{R}\left({\bf L}\Lambda^{\mathfrak{a}}(X),{\bf L}\Lambda^{\mathfrak{a}}(Y)\right) \\
 & \simeq {\bf R}\Hom_{R}\left(X,{\bf L}\Lambda^{\mathfrak{a}}(Y)\right) \\
 & \simeq {\bf R}\Hom_{R}\left({\bf R}\Gamma_{\mathfrak{a}}(X),{\bf L}\Lambda^{\mathfrak{a}}(Y)\right) \\
 & \simeq {\bf R}\Hom_{R}\left({\bf R}\Gamma_{\mathfrak{a}}(X),Y\right) \\
 & \simeq {\bf R}\Hom_{R}\left({\bf R}\Gamma_{\mathfrak{a}}(X),{\bf R}\Gamma_{\mathfrak{a}}(Y)\right).
\end{split}
\end{equation*}
\end{corollary}

\begin{prf}
By Corollary \ref{3.3.6}, Adjointness Isomorphism, and Theorem \ref{3.4.1}, we have
\begin{equation} \label{eq:1}
\begin{split}
{\bf L}\Lambda^{\mathfrak{a}}\left({\bf R}\Hom_{R}(X,Y)\right) & \simeq {\bf R}\Hom_{R}\left({\bf R}\Gamma_{\mathfrak{a}}(R),{\bf R}\Hom_{R}(X,Y)\right) \\
 & \simeq {\bf R}\Hom_{R}\left({\bf R}\Gamma_{\mathfrak{a}}(R)\otimes_{R}^{\bf L}X,Y\right) \\
 & \simeq {\bf R}\Hom_{R}\left({\bf R}\Gamma_{\mathfrak{a}}(X),Y\right) \\
 & \simeq {\bf R}\Hom_{R}\left(X,{\bf L}\Lambda^{\mathfrak{a}}(Y)\right).
\end{split}
\end{equation}
Further, by Theorem \ref{3.4.1}, \cite[Corollary on Page 6]{AJL}, and \cite[Proposition 3.2.2]{Li}, we have
\begin{equation} \label{eq:2}
\begin{split}
{\bf R}\Hom_{R}\left({\bf R}\Gamma_{\mathfrak{a}}(X),{\bf L}\Lambda^{\mathfrak{a}}(Y)\right) & \simeq {\bf R}\Hom_{R}\left({\bf R}\Gamma_{\mathfrak{a}}\left({\bf R}\Gamma_{\mathfrak{a}}(X)\right),Y\right) \\
 & \simeq {\bf R}\Hom_{R}\left({\bf R}\Gamma_{\mathfrak{a}}(X),Y\right) \\
 & \simeq {\bf R}\Hom_{R}\left({\bf R}\Gamma_{\mathfrak{a}}(X),{\bf R}\Gamma_{\mathfrak{a}}(Y)\right).
\end{split}
\end{equation}
Moreover, by Theorem \ref{3.4.1} and \cite[Corollary on Page 6]{AJL}, we have
\begin{equation} \label{eq:3}
\begin{split}
{\bf R}\Hom_{R}\left({\bf L}\Lambda^{\mathfrak{a}}(X),{\bf L}\Lambda^{\mathfrak{a}}(Y)\right) & \simeq {\bf R}\Hom_{R}\left({\bf R}\Gamma_{\mathfrak{a}}\left({\bf L}\Lambda^{\mathfrak{a}}(X)\right),Y\right) \\
 & \simeq {\bf R}\Hom_{R}\left({\bf R}\Gamma_{\mathfrak{a}}(X),Y\right).
\end{split}
\end{equation}
Combining the isomorphisms \eqref{eq:1}, \eqref{eq:2}, and \eqref{eq:3}, we get all the desired isomorphisms.
\end{prf}

Now we turn our attention to the Grothendieck's Local Duality, and demonstrate how to derive it from the Greenlees-May Duality.

We need the definition of a dualizing complex.

\begin{definition} \label{3.4.2}
A \textit{dualizing complex} \index{dualizing complex} for $R$ is an $R$-complex $D \in \mathcal{D}^{f}_{\square}(R)$ that satisfies the following conditions:
\begin{enumerate}
\item[(i)] The homothety morphism $\chi^{D}_{R}:R \rightarrow {\bf R}\Hom_{R}(D,D)$ is an isomorphism in $\mathcal{D}(R)$.
\item[(ii)] $\id_{R}(D)<\infty$.
\end{enumerate}
Moreover, if $R$ is local, then a dualizing complex $D$ is said to be \textit{normalized} \index{normalized dualizing complex} if $\sup(D)=\dim(R)$.
\end{definition}

It is clear that if $D$ is a dualizing complex for $R$, then so is $\Sigma^{s}D$ for every $s \in \mathbb{Z}$, which accounts for the non-uniqueness of dualizing complexes. Further, $\Sigma^{\dim(R)-\sup(D)}D$ is a normalized dualizing complex.

\begin{example} \label{3.4.3}
Let $(R,\mathfrak{m},k)$ be a local ring with a normalized dualizing complex $D$. Then ${\bf R}\Gamma_{\mathfrak{m}}(D) \simeq E_{R}(k)$. For a proof, refer to \cite[Proposition 6.1]{Ha}.
\end{example}

The next theorem determines precisely when a ring enjoys a dualizing complex.

\begin{theorem} \label{3.4.4}
The the following assertions are equivalent:
\begin{enumerate}
\item[(i)] $R$ has a dualizing complex.
\item[(ii)] $R$ is a homomorphic image of a Gorenstein ring of finite Krull dimension.
\end{enumerate}
\end{theorem}

\begin{prf}
See \cite[Page 299]{Ha} and \cite[Corollary 1.4]{Kw}.
\end{prf}

Now we prove the Local Duality Theorem for complexes. \index{local duality theorem for complexes}

\begin{theorem} \label{3.4.5}
Let $(R,\mathfrak{m})$ be a local ring with a dualizing complex $D$, and $X \in \mathcal{D}^{f}_{\square}(R)$. Then
$$H^{i}_{\mathfrak{m}}(X) \cong \Ext^{\dim(R)-i-\sup(D)}_{R}(X,D)^{\vee}$$
for every $i \in \mathbb{Z}$.
\end{theorem}

\begin{prf}
Clearly, we have
$$\Ext^{\dim(R)-i-\sup(D)}_{R}(X,D) \cong \Ext^{-i}_{R}\left(X,\Sigma^{\dim(R)-\sup(D)}D\right)$$
for every $i \in \mathbb{Z}$, and $\Sigma^{\dim(R)-\sup(D)}D$ is a normalized dualizing for $R$. Hence by replacing $D$ with $\Sigma^{\dim(R)-\sup(D)}D$, it suffices to assume that $D$ is a normalized dualizing complex and prove the isomorphism
$H^{i}_{\mathfrak{m}}(X) \cong \Ext^{-i}_{R}(X,D)^{\vee}$ for every $i \in \mathbb{Z}$.
By Theorem \ref{3.4.1}, we have
\begin{equation} \label{eq:3.4.5.1}
{\bf R}\Hom_{R}\left({\bf R}\Gamma_{\mathfrak{m}}(X),E_{R}(k)\right) \simeq {\bf R}\Hom_{R}\left(X,{\bf L}\Lambda^{\mathfrak{m}}\left(E_{R}(k)\right)\right).
\end{equation}
But since $E_{R}(k)$ is injective, it provides a semi-injective resolution of itself, so we have
\begin{equation} \label{eq:3.4.5.2}
{\bf R}\Hom_{R}\left({\bf R}\Gamma_{\mathfrak{m}}(X),E_{R}(k)\right) \simeq \Hom_{R}\left({\bf R}\Gamma_{\mathfrak{m}}(X),E_{R}(k)\right).
\end{equation}
Besides, by Example \ref{3.4.3}, \cite[Corollary on Page 6]{AJL}, and \cite[Proposition 2.7]{Fr}, we have
\begin{equation} \label{eq:3.4.5.3}
\begin{split}
{\bf L}\Lambda^{\mathfrak{m}}\left(E_{R}(k)\right) & \simeq {\bf L}\Lambda^{\mathfrak{m}}\left({\bf R}\Gamma_{\mathfrak{m}}(D)\right) \\
 & \simeq {\bf L}\Lambda^{\mathfrak{m}}(D) \\
 & \simeq D\otimes_{R}^{\bf L}\widehat{R}^{\mathfrak{m}} \\
 & \simeq D\otimes_{R}\widehat{R}^{\mathfrak{m}}.
\end{split}
\end{equation}
Combining \eqref{eq:3.4.5.1}, \eqref{eq:3.4.5.2}, and \eqref{eq:3.4.5.3}, we get
$$\Hom_{R}\left({\bf R}\Gamma_{\mathfrak{m}}(X),E_{R}(k)\right) \simeq {\bf R}\Hom_{R}\left(X,D\otimes_{R}\widehat{R}^{\mathfrak{m}}\right).$$
Taking Homology, we obtain
\begin{equation} \label{eq:3.4.5.4}
\begin{split}
\Hom_{R}\left(H^{i}_{\mathfrak{m}}(X),E_{R}(k)\right) & \cong \Hom_{R}\left(H_{-i}\left({\bf R}\Gamma_{\mathfrak{m}}(X)\right),E_{R}(k)\right) \\
 & \cong H_{i}\left(\Hom_{R}\left({\bf R}\Gamma_{\mathfrak{m}}(X),E_{R}(k)\right)\right) \\
 & \cong H_{i}\left({\bf R}\Hom_{R}\left(X,D\otimes_{R}\widehat{R}^{\mathfrak{m}}\right)\right) \\
 & \cong \Ext^{-i}_{R}\left(X,D\otimes_{R}\widehat{R}^{\mathfrak{m}}\right) \\
\end{split}
\end{equation}
for every $i \in \mathbb{Z}$.

Since $X \in \mathcal{D}^{f}_{\square}(R)$, we have $X\otimes_{R}\widehat{R}^{\mathfrak{m}} \in \mathcal{D}^{f}_{\square}\left(\widehat{R}^{\mathfrak{m}}\right)$, so $H^{i}_{\mathfrak{m}\widehat{R}^{\mathfrak{m}}}\left(X\otimes_{R}\widehat{R}^{\mathfrak{m}}\right)$ is an artinian $\widehat{R}^{\mathfrak{m}}$-module by \cite[Proposition 2.1]{HD}, and thus Matlis reflexive for every $i \in \mathbb{Z}$. Moreover, $D\otimes_{R}\widehat{R}^{\mathfrak{m}}$ is a normalized dualizing complex for $\widehat{R}^{\mathfrak{m}}$. Therefore, using the isomorphism \eqref{eq:3.4.5.4} over the $\mathfrak{m}$-adically complete ring $\widehat{R}^{\mathfrak{m}}$, we obtain
\begin{equation} \label{eq:3.4.5.5}
\begin{split}
H_{\mathfrak{m}}^{i}(X) & \cong H_{\mathfrak{m}}^{i}(X)\otimes_{R} \widehat{R}^{\mathfrak{m}}\\
 & \cong H^{i}_{\mathfrak{m}\widehat{R}^{\mathfrak{m}}}\left(X\otimes_{R} \widehat{R}^{\mathfrak{m}}\right) \\
 & \cong \Hom_{\widehat{R}^{\mathfrak{m}}}\left(\Hom_{\widehat{R}^{\mathfrak{m}}}\left(H^{i}_{\mathfrak{m}\widehat{R}^{\mathfrak{m}}}\left(X\otimes_{R} \widehat{R}^{\mathfrak{m}}\right),E_{\widehat{R}^{\mathfrak{m}}}(k)\right),E_{\widehat{R}^{\mathfrak{m}}}(k)\right) \\
 & \cong \Hom_{\widehat{R}^{\mathfrak{m}}}\left(\Ext^{-i}_{\widehat{R}^{\mathfrak{m}}}\left(X\otimes_{R}\widehat{R}^{\mathfrak{m}},D\otimes_{R}\widehat{R}^{\mathfrak{m}} \right),E_{\widehat{R}^{\mathfrak{m}}}(k)\right) \\
 & \cong \Hom_{\widehat{R}^{\mathfrak{m}}}\left(\Ext^{-i}_{R}(X,D)\otimes_{R}\widehat{R}^{\mathfrak{m}},E_{\widehat{R}^{\mathfrak{m}}}(k)\right) \\
\end{split}
\end{equation}
for every $i \in \mathbb{Z}$. However, ${\bf R}\Hom_{R}(X,D) \in \mathcal{D}^{f}_{\sqsubset}(R)$, so $\Ext^{-i}_{R}(X,D)$ is a finitely generated $R$-module for every $i \in \mathbb{Z}$. It follows that
\begin{equation} \label{eq:3.4.5.6}
\Hom_{\widehat{R}^{\mathfrak{m}}}\left(\Ext^{-i}_{R}(X,D)\otimes_{R}\widehat{R}^{\mathfrak{m}},E_{\widehat{R}^{\mathfrak{m}}}(k)\right) \cong \Hom_{R}\left(\Ext^{-i}_{R}(X,D),E_{R}(k)\right)
\end{equation}
for every $i \in \mathbb{Z}$.
Combining \eqref{eq:3.4.5.5} and \eqref{eq:3.4.5.6}, we obtain
$$H^{i}_{\mathfrak{m}}(X) \cong \Hom_{R}\left(\Ext^{-i}_{R}(X,D),E_{R}(k)\right)$$
for every $i \in \mathbb{Z}$ as desired.
\end{prf}

Our next goal is to obtain the Local Duality Theorem for modules. But first we need the definition of a dualizing module.

\begin{definition} \label{3.4.6}
Let $(R,\mathfrak{m})$ be a local ring. A \textit{dualizing module} \index{dualizing module} for $R$ is a finitely generated $R$-module $\omega$ that satisfies the following conditions:
\begin{enumerate}
\item[(i)] The homothety map $\chi^{\omega}_{R}:R \rightarrow \Hom_{R}(\omega,\omega)$, given by $\chi^{\omega}_{R}(a)=a1^{\omega}$ for every $a \in R$, is an isomorphism.
\item[(ii)] $\Ext^{i}_{R}(\omega,\omega)=0$ for every $i \geq 1$.
\item[(iii)] $\id_{R}(\omega)<\infty$.
\end{enumerate}
\end{definition}

The next theorem determines precisely when a ring enjoys a dualizing module.

\begin{theorem} \label{3.4.7}
Let $(R,\mathfrak{m})$ be a local ring. Then the following assertions are equivalent:
\begin{enumerate}
\item[(i)] $R$ has a dualizing module.
\item[(ii)] $R$ is a Cohen-Macaulay local ring which is a homomorphic image of a Gorenstein local ring.
\end{enumerate}
Moreover in this case, the dualizing module is unique up to isomorphism.
\end{theorem}

\begin{prf}
See \cite[Corollary 2.2.13]{Wa} and \cite[Theorem 3.3.6]{BH}.
\end{prf}

Since the dualizing module for $R$ is unique whenever it exists, we denote a choice of the dualizing module by $\omega_{R}$.

\begin{proposition} \label{3.4.8}
Let $(R,\mathfrak{m})$ be a Cohen-Macaulay local ring, and $\omega$ a finitely generated $R$-module. Then the following assertions are equivalent:
\begin{enumerate}
\item[(i)] $\omega$ is a dualizing module for $R$.
\item[(ii)] $\omega^{\vee} \cong H^{\dim(R)}_{\mathfrak{m}}(R)$.
\end{enumerate}
\end{proposition}

\begin{prf}
See \cite[Definition 12.1.2, Exercises 12.1.23 and 12.1.25, and Remark 12.1.26]{BS}, and \cite[Definition 3.3.1]{BH}.
\end{prf}

We can now derive the Local Duality Theorem for modules. \index{local duality theorem for modules}

\begin{theorem} \label{3.4.9}
Let $(R,\mathfrak{m})$ be a local ring with a dualizing module $\omega_{R}$, and $M$ a finitely generated $R$-module. Then
$$H^{i}_{\mathfrak{m}}(M) \cong \Ext^{\dim(R)-i}_{R}(M,\omega_{R})^{\vee}$$
for every $i \geq 0$.
\end{theorem}

\begin{prf}
By Theorem \ref{3.4.7}, $R$ is a Cohen-Macaulay local ring which is a homomorphic image of a Gorenstein local ring $S$. Since $S$ is local, we have $\dim(S)<\infty$. Hence Theorem \ref{3.4.4} implies that $R$ has a dualizing complex $D$. Since $R$ is Cohen-Macaulay, we have $H^{i}_{\mathfrak{m}}(R)=0$ for every $i \neq \dim(R)$. On the other hand, by Theorem \ref{3.4.5}, we have
\begin{equation} \label{eq:3.4.9.1}
\begin{split}
H^{i}_{\mathfrak{m}}(R) & \cong \Ext^{\dim(R)-i-\sup(D)}_{R}(R,D)^{\vee} \\
 & \cong H_{-\dim(R)+i+\sup(D)}\left({\bf R}\Hom_{R}(R,D)\right)^{\vee} \\
 & \cong H_{-\dim(R)+i+\sup(D)}(D)^{\vee}.
\end{split}
\end{equation}
It follows from the display \eqref{eq:3.4.9.1} that $H_{-\dim(R)+i+\sup(D)}(D)=0$ for every $i \neq \dim(R)$, i.e. $H_{i}(D)=0$ for every $i \neq \sup(D)$. Therefore, we have $D \simeq \Sigma^{\sup(D)}H_{\sup(D)}(D)$. In addition, letting $i=\dim(R)$ in the display \eqref{eq:3.4.9.1}, we get $H^{\dim(R)}_{\mathfrak{m}}(R) \cong H_{\sup(D)}(D)^{\vee}$, which implies that $\omega_{R} \cong H_{\sup(D)}(D)$ by Proposition \ref{3.4.8}. It follows that $D \simeq \Sigma^{\sup(D)}\omega_{R}$.

Now let $M$ be a finitely generated $R$-module. Then by Theorem \ref{3.4.5}, we have
\begin{equation*}
\begin{split}
H^{i}_{\mathfrak{m}}(M) & \cong \Ext^{\dim(R)-i-\sup(D)}_{R}(M,D)^{\vee} \\
 & \cong H_{-\dim(R)+i+\sup(D)}\left({\bf R}\Hom_{R}\left(M,D\right)\right)^{\vee} \\
 & \cong H_{-\dim(R)+i+\sup(D)}\left({\bf R}\Hom_{R}\left(M,\Sigma^{\sup(D)}\omega_{R}\right)\right)^{\vee} \\
 & \cong H_{-\dim(R)+i}\left({\bf R}\Hom_{R}\left(M,\omega_{R}\right)\right)^{\vee} \\
 & \cong \Ext^{\dim(R)-i}_{R}\left(M,\omega_{R}\right)^{\vee}. \\
\end{split}
\end{equation*}
\end{prf}


\end{document}